\documentclass{elsarticle}

\usepackage{amsmath,amssymb,amsfonts}
\usepackage{algorithmic}
\usepackage{graphicx}
\usepackage{textcomp}
\usepackage{xcolor}
\usepackage{enumitem}
\usepackage{xcolor}
\usepackage[normalem]{ulem}
\usepackage{steinmetz}
\usepackage{multirow}
\usepackage{tikz}
\usepackage{caption}
\usepackage{subcaption}
\usepackage{url}

\newcommand{\R}{\color{black}}

\begin{document}
	

\title{Power Distribution Network Reconfiguration for Distributed Generation Maximization}

\author[1]{Kin Cheong Sou\corref{cor1} \fnref{fn1}} \ead{sou12@mail.nsysu.edu.tw}
\author[2]{Gabriel Malmer} \ead{gabriel.malmer@iea.lth.se}
\author[2]{Lovisa Thorin} \ead{lovisathorin@gmail.com}
\author[2]{Olof Samuelsson} \ead{olof.samuelsson@iea.lth.se}

\cortext[cor1]{Corresponding author}
\fntext[fn1]{K.C.~Sou is partially supported by the National Science and Technology Council (NSTC) of Taiwan: 112-2221-E-110-007-, 113-2221-E-110-022-MY2 and 114-2218-E-007-011-}

\affiliation[1]{organization={National Sun Yat-sen University}, addressline={No.~70, Lien-Hai Rd.}, city={Kaohsiung}, country={Taiwan}}
	
\affiliation[2]{organization={Lund University}, addressline={Box 118},
		city={Lund},
		postcode={221 00},
		country={Sweden}}

\begin{abstract}
    Network reconfiguration can significantly increase the hosting capacity (HC) for distributed generation (DG) in radially operated systems, thereby reducing the need for costly infrastructure upgrades. However, when the objective is DG maximization, jointly optimizing topology and power dispatch remains computationally challenging. Existing approaches often rely on relaxations or approximations, yet we provide counterexamples showing that interior point methods, linearized DistFlow and second-order cone relaxations all yield erroneous results. To overcome this, we propose a solution framework based on the exact DistFlow equations, formulated as a bilinear program and solved using spatial branch-and-bound (SBB). Numerical studies on standard benchmarks and a 533-bus real-world system demonstrate that our proposed method reliably performs reconfiguration and dispatch within time frames compatible with real-time operation.
\end{abstract}

\maketitle

\section{Introduction}
Distribution networks are generally built in a looped fashion, but operated with one point normally-open in each loop resulting in radial networks. The set of switches to leave open (hence the network configuration) offers a degree-of-freedom that is normally not exploited in the industry. In this paper, we study the joint optimization of network reconfiguration and power control of distributed generation (DG), and we call this the \emph{reconfiguration optimal power flow} (ROPF) problem, generalizing optimal power flow (OPF). While enforcing network configuration and operational constraints (e.g., voltage and line current limits and DG restrictions), the ROPF problem traditionally seeks to minimize objectives including total line loss, fuel cost or voltage deviation (e.g., \cite{baran1989network,jabr2012minimum,taylor2012convex,ababei2010efficient,santos2022dynamic,bahrami2024dynamic}). However, in this paper we consider the maximization of total DG output in ROPF as an emerging topic for renewable energy integration (e.g., \cite{capitanescu2014assessing,fu2018toward,song2019new,home2022increasing}). The DG maximization objective has a crucial effect on the tractability of the ROPF problem and the usefulness of the resulting dispatch solution. {\R This is the first main discovery of this paper.}

Diverse solution methods for the ROPF problem have been explored. For instance, heuristics and rule-based approaches have been considered due to their ease for online computation (e.g., \cite{baran1989network,ababei2010efficient,hedman2011review}). To more explicitly search for the optimal solution, the ROPF problem is formulated and solved as an optimization problem, akin to OPF's solution approach. A variety of methods have been proposed to solve the ROPF optimization problem. Meta-heuristics algorithms (e.g., \cite{su2003network,zhang2007improved,fu2018toward}) are commonly considered due to their flexibility and ease of use. However, lack of optimality guarantee and unpredictable runtime are major drawbacks of these methods. Thus, more direct optimization methods are explored. For example, using DC power flow equations or linearized DistFlow equations to approximate the power system physics, \cite{ramos2005path,fisher2008optimal,ergun2012transmission,singh2022joint} investigate mixed integer linear programming (MILP) models for ROPF. Alternatively, to directly include the AC power flow equations, mixed integer conic programming (MICP) ROPF models are proposed (e.g., \cite{jabr2012minimum,farivar2013branch,taylor2012convex,kocuk2017new}). Theoretical support and numerical demonstrations have been reported (e.g., \cite{farivar2013branch,low2014convex,jabr2012minimum,taylor2012convex,lee2014robust,home2022increasing}). However, it is not well recognized that the objective of the ROPF problem plays an important role in the suceess of the solution approaches. In particular, except for \cite{home2022increasing} all aformentioned previous works consider ``classical'' objectives such as loss minimization, convex quadratic fuel cost minimization, conservation voltage reduction, etc. For the ``non-classical'' objective of DG maximization in this paper, it turns out that popular methods such as the MICP apporach and interior point algorithm may lead to undesirable control dispatch violating grid constraints such as voltage and current limits. We support our claim with a repeatable case study involving a three-bus example with all computation details presented. To the best of our knowledge, the adverse effects of these pitfalls have not been well recognized in the literature.

The difficulty with DG maximizing ROPF necessitates a solution method with acceptable accuracy and efficiency. With the only nonlinearity in the model being polynomials (due to the power flow equations), the ROPF problem can be interpreted as a bilinear programming problem using specialized solver (e.g., Gurobi). Similar approaches have been studied in \cite{lavorato2011imposing,capitanescu2014assessing}, where the ROPF problem is treated as a more general mixed integer nonlinear programming (MINLP) problem and solved using the spatial branch-and-bound (SBB) algorithm. In terms of accuracy, the SBB algorithm is sufficient since it guarantees optimality regardless of the ROPF objective type. However, since SBB is fundamentally an exhaustive search it is important to investigate its computational efficiency for the DG maximizing ROPF problem. The rather long runtime reported in \cite{lavorato2011imposing,capitanescu2014assessing} does not imply efficiency especially for real-time applications. The efficiency of SBB is significantly affected by the choice of problem formulation and algorithm implementation. A crucial factor in problem formulation is the model of power system physics described by the power flow equations. Typical forms of the power flow equations include, for example, AC power flow equations with voltages in polar coordinates (e.g., \cite{lavorato2011imposing}) and voltages in rectangular coordinates (e.g., \cite{capitanescu2014assessing}), branch flow equations \cite{farivar2013branch}, DistFlow equations \cite{baran1989network}, modified AC power flow equations (e.g., \cite{exposito1999reliable,jabr2006radial}), etc. On the other hand, different SBB implementations (e.g., Gurobi, BARON, Knitro) behave differently because they make different assumptions about the problem structure. Through our case studies, we identify the combination of DistFlow equations and Gurobi (version 9 or later, with recently added capability to handle bilinear programs) as a promising choice of solution method for the DG maximizing ROPF problem. For instance, for a 533-bus real example in Sweden, ROPF can be performed within 10 minutes when allowing up to four switch status changes. To the best of our knowledge, successful demonstration of SBB to solve ROPF models with the exact DistFlow equations is not well-known, even though the DistFlow equations are commonly relaxed and used in the MICP approach (e.g., \cite{farivar2013branch}).  We emphasize that the appropriate choice of formulation and SBB solver is important. We will provide counterexamples to illustrate the poor computational performance of some improper choices. {\R We will also demonstrate that the DistFlow equations play a vital role in the efficient computation of the ROPF problem. This is the second main discovery of this paper.}

    {\R The main novelties and contributions of this paper are listed as follows:}
	\begin{itemize}
		\item We raise the awareness that the DG maximizing objective fundamentally changes the tractability of the ROPF problem and illustrate the pitfalls of possible overextension of popular solution methods such as MICP relaxation, interior point algorithms and linearization approximation.
		
		\item We point out the importance of the choice of problem formulation and optimization solver for the DG maximizing ROPF problem, by illustrating a promising choice and other improper choices.
	\end{itemize}

The rest of this paper is organized as follows: Section~\ref{sec:problem_description} introduces the notation, defines the ROPF optimization model and describes the proposed solution approach. Section~\ref{sec:demo} illustrates and evaluates the performance of the proposed ROPF approach. Section~\ref{sec:pitfalls} justifies the proposed method by showing pitfalls encountered by common alternatives for DG maximizing ROPF. Section~\ref{sec:real_network} discusses a practical renewable energy hosting capacity maximization problem for a 533-bus distribution network in Sweden using the proposed method.

\section{Problem Description} \label{sec:problem_description}

\subsection{Problem Overview}
We consider balanced power distribution systems with line shunt elements ignored so that the system can be described by the DistFlow equations \cite{baran1989network}. The lines are equipped with switches which can be closed or open. The network is built with loop(s) but is operated in a radial topology. Power exchange is possible at the slack bus. In addition, load and DG are present in the system. The loads are assumed to be known. The DG sources are interfaced with inverters where both active and reactive power outputs are considered adjustable. By jointly adjusting network configuration and inverter outputs, the ROPF problem considered in this paper seeks to \emph{maximize} the total distributed active power generation in the system with the following constraints: (a) DistFlow equations, (b) network radiality and allowable switch status changes, (c) inverter power output bounds, apparent power capacity, power factor and (d) grid constraints (voltage and line current limits).

We comment that our adopted DistFlow model, while restricted by the three-phase balance assumption, is well accepted as evidenced by numerous influential results built upon the same assumption (e.g., \cite{baran1989network,jabr2012minimum,taylor2012convex}). Also, it is common to find representative distribution system models described exactly by the DistFlow equations (e.g., MATPOWER \cite{zimmerman2010matpower}, REDS repository \cite{ababei2010efficient}, EU examples \cite{mateo2018european}).

%
%
%


\subsection{Notation and Assumptions} \label{subsec:notations}
\begin{itemize}
	\item Bus 0 is the slack bus with voltage $v_0 = 1 \phase{0^\circ}$ pu. The remaining buses are labelled $1, 2, \ldots, N$. We denote $\mathcal{N} := \{1, \ldots, N \}$ and $\bar{\mathcal{N}} := \mathcal{N} \cup \{0\}$.
	
	\item For any bus $n \in \bar{\mathcal{N}}$, the voltage is $v_n$. We define $\nu_n := |v_n|^2$ as the squared voltage for bus $n$. $\nu_n$ is a state decision variable in the ROPF problem.
	
	\item For any bus $n \in \mathcal{N}$, the active and reactive power loads are $p^L_n$ and $q^L_n$. Their values are assumed to be known.
	
	\item For any bus $n \in \mathcal{N}$, the active and reactive power generation are $p^G_n$ and $q^G_n$. They are control decision variables of the ROPF problem.
	
	\item In a radial network, each bus $n \in \mathcal{N}$ is associated with exactly one line, denoted by line $n$, with reference direction pointing away from bus 0. The current and power flow follow the same reference direction of the line.
	
	\item $\pi_n \in \bar{\mathcal{N}}$ is the ``parent'' bus of $n \in \mathcal{N}$ (so that line $n$ goes from $\pi_n$ to $n$).
	
	\item For any $n \in \bar{\mathcal{N}}$, $c(n) \subseteq \mathcal{N}$ is the set of buses with parent being $n$. That is, $c(n) = \{k \in \mathcal{N} \mid \pi_k = n\}$.
	
	\item For any bus $n \in \mathcal{N}$, $P_n$ is the active power flow from bus $\pi_n$ to bus $n$ measured at bus $\pi_n$. $Q_n$ is defined similarly for reactive power flow. $P_n$ and $Q_n$ are state decision variables of the ROPF problem.
	
	\item Line shunt elements are ignored. For any bus $n \in \mathcal{N}$, $\ell_n$ is the squared magnitude current from bus $\pi_n$ to bus $n$ (also the value from $n$ to $\pi_n$). $\ell_n$ is a state decision variable in the ROPF problem.
	
	
	\item For a line $\{n,k\}$ in a reconfigurable network, the series resistance and reactance are $r_{nk}$ and $x_{nk}$ respectively. Also, the square line current upper limit is denoted by $\overline{\ell}_{nk}$. Note that $r_{nk} = r_{kn}$, $x_{nk} = x_{kn}$ and $\overline{\ell}_{nk} = \overline{\ell}_{kn}$.
	
	\item For a reconfigurable network, $\pi_{nk} \in \{0,1\}$ with $(n,k) \in \bar{\mathcal{N}} \times \bar{\mathcal{N}}$ such that $\pi_{nk} = 1$ if and only if bus $n$ is the parent of bus $k$ (i.e., $\pi_k = n$). $\pi_{nk}$ is a ROPF control decision variable. Also, we define the auxiliary decision variable $\alpha_{nk} = \pi_{nk} + \pi_{kn}$ for $(n,k) \in \bar{\mathcal{N}} \times \bar{\mathcal{N}}$ to encode the switch status of the line between bus $n$ and bus $k$.
	
\end{itemize}

\subsection{ROPF Problem Constraints Derivation}

\subsubsection{DistFlow Equations for Fixed Network Configuration}

With fixed network configuration, the physics of distribution systems is described by the DistFlow equations \cite{baran1989network}:
\begin{subequations} \label{eqn:DistFlow}
	\begin{align}
		p^G_n - p^L_n &= -P_n + \sum\limits_{k \in c(n)} P_k + r_n \ell_n,  & n \in \mathcal{N} \label{eqn:DistFlow_P} \\
		q^G_n - q^L_n &= -Q_n + \sum\limits_{k \in c(n)} Q_k + x_n \ell_n, & n \in \mathcal{N} \label{eqn:DistFlow_Q} \\
		\nu_{\pi_n} - \nu_n &= 2 r_n P_n + 2 x_n Q_n - (r_n^2 + x_n^2) \ell_n, & n \in \mathcal{N} \label{eqn:DistFlow_Ohm} \\
		\nu_{\pi_n} \ell_n &= P_n^2 + Q_n^2, & n \in \mathcal{N} \label{eqn:DistFlow_SVI}
	\end{align}
\end{subequations}
where $r_n$ and $x_n$ are the series resistance and series reactance of line $n$, respectively. The set $c(n)$ is defined in Section~\ref{subsec:notations}.


\subsubsection{Network Configuration Modelling} \label{subsec:network_configuration} We impose bounds on the switch status decision variables $\alpha_{nk}$ defined in Section~\ref{subsec:notations} to indicate always closed and always open switches:
\begin{equation} \label{eqn:alpha_bounds}
	\underline{\alpha}_{nk} \le \alpha_{nk} \le \overline{\alpha}_{nk}, \qquad (n,k) \in \bar{\mathcal{N}} \times \bar{\mathcal{N}}
\end{equation}
A network with $N+1$ buses is radial if and only if it is connected and contains exactly $N$ lines. Radiality can be imposed by the (necessary and sufficient) single-commodity flow conditions. Using the parent decision variables $\pi_{nk}$ defined in Section~\ref{subsec:notations} and auxiliary continuous ``flow'' variables $f_{ij}$ for $(i,j) \in \bar{\mathcal{N}} \times \bar{\mathcal{N}}$, these conditions can be described as \cite{sou2022joint,sou2022relaxed}
\begin{subequations} \label{eqn:SCF}
	\begin{align}
		& \sum\limits_{k  \in \bar{\mathcal{N}}} f_{nk} - \sum\limits_{k  \in \bar{\mathcal{N}}} f_{kn} = -1, & n \in \mathcal{N} \label{eqn:flow_conservation} \\
		& 0 \le f_{nk} \le N \pi_{nk}, & (n,k) \in \bar{\mathcal{N}} \times \bar{\mathcal{N}} \label{eqn:flow_enabling} \\
		& \sum\limits_{(n,k) \in \bar{\mathcal{N}} \times \bar{\mathcal{N}}} \pi_{nk} = N. \label{eqn:tree_num_of_edges}
	\end{align}
\end{subequations}

We may wish to restrict the number of switch status changes denoted by $K$. Let parameters $\alpha^0_{nk} \in \{0,1\}$ for $(n,k) \in \bar{\mathcal{N}} \times \bar{\mathcal{N}}$ denote the initial network configuration such that $\alpha_{nk} = \alpha_{kn} = 1$ if switch $\{n,k\}$ is initially closed and $\alpha_{nk} = \alpha_{kn} = 0$ if $\{n,k\}$ is initially open. Then, the desired constraint is
\begin{equation} \label{eqn:status_change}
	\sum\limits_{(n,k) \in \bar{\mathcal{N}} \times \bar{\mathcal{N}}} \!\! \Big( \alpha^0_{nk} (1 - \alpha_{nk}) + \alpha_{nk} (1 - \alpha^0_{nk}) \Big) \le 2 K.
\end{equation}
Note that when a switch is closed another switch must be opened simultaneously to keep the radial network configuration. Hence, $K$ is typically an even number.

\subsubsection{DistFlow Equations for Configurable Network} \label{subsec:reconfig_DistFlow}
When the network configuration is parameterized by decision variable $\pi$ (also $\alpha$), certain configuration related terms in \eqref{eqn:DistFlow} should be generalized. For instance, the children set $c(n)$ in \eqref{eqn:DistFlow_P} and \eqref{eqn:DistFlow_Q} is not determined before the problem is solved. Thus, the following changes are needed in \eqref{eqn:DistFlow_P} and \eqref{eqn:DistFlow_Q}
\begin{displaymath}
	\sum\limits_{k \in c(n)} P_k \rightarrow \sum\limits_{k \in \mathcal{N}} \pi_{nk} P_k, \;\; \sum\limits_{k \in c(n)} Q_k \rightarrow \sum\limits_{k \in \mathcal{N}} \pi_{nk} Q_k, \;\; n \in \mathcal{N}
\end{displaymath}
to describe the sum of power flows leaving bus $n \in \mathcal{N}$. Similarly, while line $n$ always points to bus $n$, the parent $\pi_n$ is not fixed before optimization. The line parameters and parent bus voltage in \eqref{eqn:DistFlow_P} to \eqref{eqn:DistFlow_SVI} become
\begin{displaymath}
	r_n \rightarrow \sum\limits_{k \in \bar{\mathcal{N}}} \pi_{kn} r_{kn}, \;\; x_n \rightarrow \sum\limits_{k \in \bar{\mathcal{N}}} \pi_{kn} x_{kn}, \;\; \nu_{\pi_n} \rightarrow \sum\limits_{k \in \bar{\mathcal{N}}} \pi_{kn} \nu_k,
\end{displaymath}
where $r_{kn}$ and $x_{kn}$ are the series resistance and reactance of line $\{n,k\}$ as defined in Section~\ref{subsec:notations}. Hence, the DistFlow equations for configurable networks become
\begin{subequations} \label{eqn:DistFlow_reconfig}
	\begin{align}
		& p^G_n \! - \! p^L_n = \! -P_n \!\! + \!\! \sum\limits_{k \in \mathcal{N}} \pi_{nk} P_k \! + \!\! \Big(\sum\limits_{k \in \bar{\mathcal{N}}} \pi_{kn} r_{kn}\Big) \ell_n, \; n \!\in\! \mathcal{N} \label{eqn:DistFlow_reconfig_P} \\
		& q^G_n \!- \! q^L_n = \! -Q_n \!\!+\!\! \sum\limits_{k \in \mathcal{N}} \pi_{nk} Q_k \!+\!\! \Big(\sum\limits_{k \in \bar{\mathcal{N}}} \!\pi_{kn} x_{kn}\Big) \ell_n,  n \!\in\! \mathcal{N} \label{eqn:DistFlow_reconfig_Q} \\
		& \sum\limits_{k \in \bar{\mathcal{N}}} \pi_{kn} \nu_k - \nu_n = 2 (\sum\limits_{k \in \bar{\mathcal{N}}} \pi_{kn} r_{kn}) P_n + 2 (\sum\limits_{k \in \bar{\mathcal{N}}} \pi_{kn} x_{kn}) Q_n \nonumber \\
		& \quad - \Big(\sum\limits_{k \in \bar{\mathcal{N}}} \pi_{kn} (r_{kn}^2 + x_{kn}^2)\Big) \ell_n, \qquad n \in \mathcal{N}, \label{eqn:DistFlow_reconfig_Ohm} \\
		& \Big(\sum\limits_{k \in \bar{\mathcal{N}}} \pi_{kn} \nu_k \Big) \ell_n = P_n^2 + Q_n^2, \qquad n \in \mathcal{N}. \label{eqn:DistFlow_reconfig_SVI}
	\end{align}
\end{subequations}
Contrary to the fixed configuration setting in \eqref{eqn:DistFlow}, the line variables $P_n$, $Q_n$ and $\ell_n$ in \eqref{eqn:DistFlow_reconfig} do not correspond to a predefined line in the network. These line variables can be associated with any line incident to bus $n$. The actual association is determined by the values of the $\pi$ variables at optimality. {\R We note that the formulation in \eqref{eqn:DistFlow_reconfig} is different from standard reconfiguration-type formulation in the literature, where big-M type status changes are directly imposed on \eqref{eqn:DistFlow_reconfig_Ohm}.}

In \eqref{eqn:DistFlow_reconfig_P}, the bilinear term $\pi_{nk} P_k$ with 0-1 binary variable $\pi_{nk}$ and continuous variable $P_k$ can be equivalently replaced by an auxiliary variable with additional linear constraints as follows. Suppose the continuous variable $P_k$ is bounded such that $\underline{P}_k \le P_k \le \overline{P}_k$. With auxiliary continuous variable $z$ satisfying
\begin{subequations} \label{eqn:bi2lin}
	\begin{align}
		& \underline{P}_k \pi_{nk} \le z \le \overline{P}_k \pi_{nk}, \label{eqn:bi2lin_1} \\
		& P_k -\overline{P}_k (1 - \pi_{nk}) \le z \le P_k - \underline{P}_k (1- \pi_{nk}), \label{eqn:bi2lin_2}
	\end{align}
\end{subequations}
the bilinear term $\pi_{nk} P_k$ can be replaced by $z$. The claim can be verified by separately considering the cases when $\pi_{nk} = 0$ and $\pi_{nk} = 1$. 
Other binary-continuous bilinear terms $\pi_{kn} \ell_n$, $\pi_{nk}Q_k$ and $\pi_{kn} \nu_k$ in \eqref{eqn:DistFlow_reconfig} can be replaced by auxiliary variables in a similar fashion. However, in \eqref{eqn:DistFlow_reconfig_SVI} the terms $(\pi_{kn} \nu_k) \ell_n$ remain bilinear even after the replacement of $\pi_{kn} \nu_k$. This complicates the solution process of the ROPF problem. 

\subsubsection{Inverter Operational Constraints} \label{subsec:inverter}
The inverter operational constraints include (a) generation bounds, (b) apparent power rating, and (c) minimum power factor. These constraints are
\begin{subequations} \label{eqn:inverter}
	\begin{align}
		& \underline{p}^G_n \le p^G_n \le \overline{p}^G_n, \quad \underline{q}^G_n \le q^G_n \le \overline{q}^G_n, & n \in \mathcal{N} \label{eqn:inverter_pg_bounds} \\
		& (p^G_n)^2 + (q^G_n)^2 \le (C_n)^2, & n \in \mathcal{N} \label{eqn:inverter_capacity} \\
		& |q^G_n| \le \tan(\cos^{-1}(\text{pf}_n)) p^G_n, \label{eqn:inverter_pf} & n \in \mathcal{N}
	\end{align}
\end{subequations}
where $\underline{p}^G_n$, $\overline{p}^G_n$, $\underline{q}^G_n$, $\overline{q}^G_n$, $C_n$, $\text{pf}_n$ are given parameters for active and reactive power generation bounds, inverter apparent power rating and minimum power factor for bus $n$, respectively. 


\subsubsection{Grid Constraints} \label{subsec:grid_code}
For operational safety, voltage limits are imposed:
\begin{equation} \label{eqn:v_limits}
	\underline{\nu}_n \le \nu_n \le \overline{\nu}_n, \qquad n \in \bar{\mathcal{N}}
\end{equation}
where $\underline{\nu}_n$ and $\overline{\nu}_n$ are given. Common values are $\sqrt{\underline{\nu}_n} = 0.95$ pu, $\sqrt{\overline{\nu}_n} = 1.05$ pu. Also, line current limit constraints are
\begin{equation} \label{eqn:l_limit_heterogeneous}
	0 \le \ell_n \le \sum\limits_{k \in \bar{\mathcal{N}}} \pi_{kn} \overline{\ell}_{kn}, \qquad n \in \mathcal{N}
\end{equation}
where $\overline{\ell}_{kn}$ is the given square current upper limit for line $\{k,n\}$. Typically, $\sqrt{\overline{\ell}_{kn}}$ is up to a few hundred amperes.

The replacement of the bilinear terms $\pi_{nk} P_k$ in \eqref{eqn:DistFlow_reconfig_P} and $\pi_{nk} Q_k$ in \eqref{eqn:DistFlow_reconfig_Q} requires variable bounds on the line power flows $P_k$ and $Q_k$. If these bounds are not provided, they can be derived as follows:
\begin{align}
	& - \sqrt{\nu^U_k {\ell}^U_k} \le P_k \le \sqrt{\nu^U_k {\ell}^U_k}, & k \in \mathcal{N} \label{eqn:Pline_UB} \\
	& - \sqrt{\nu^U_k {\ell}^U_k} \le Q_k \le \sqrt{\nu^U_k {\ell}^U_k}, & k \in \mathcal{N} \label{eqn:Qline_UB}
\end{align}
where $\nu^U_k := \max\limits_{(n,k) \mid \overline{\alpha}_{nk} = 1} \overline{\nu}_n$ and ${\ell}^U_k := \max\limits_{(n,k) \mid \overline{\alpha}_{nk} = 1} \overline{\ell}_{nk}$. The bounds in \eqref{eqn:Pline_UB} and \eqref{eqn:Qline_UB} follow necessarily from \eqref{eqn:v_limits} and \eqref{eqn:l_limit_heterogeneous}. $P_k + j Q_k = v_{\pi_k} (I_k)^*$ with $v$ and $I$ being the complex voltage and line current. Thus, $|P_k| \le |P_k + j Q_k| \le |v_{\pi_k}| |I_k| \le \sqrt{\nu^U_k {\ell}^U_k}$. This establishes \eqref{eqn:Pline_UB}. The same argument shows \eqref{eqn:Qline_UB}. 

\subsection{DG maximizing ROPF Model and the Proposed Solution Approach} \label{subsec:problem_summary}
Maximizing the total DG active power output, the ROPF problem studied in this paper can be summarized as
\begin{align} \label{opt:reconfig}
	\begin{split}
		\underset{\substack{\nu, \ell, p^G, q^G,  P, Q, \pi, \alpha}}{\text{maximize}} & \quad \sum\limits_{n \in \mathcal{N}} p^G_n \\
		\text{subject to} & \quad \text{\eqref{eqn:alpha_bounds} to \eqref{eqn:DistFlow_reconfig}, \eqref{eqn:inverter} to \eqref{eqn:l_limit_heterogeneous}}
	\end{split}
\end{align}
In \eqref{opt:reconfig}, the decision variables without subscripts denote the vectors or matrices stacked by the corresponding variables for individual buses or lines (e.g., $p^G$ is the $N$-vector of all $p^G_n$ for $n \in \mathcal{N}$). $\nu, \ell, p^G, q^G, P, Q$ are continuous vector variables, while $\pi$ and $\alpha$ are 0-1 binary matrix variables. Except for \eqref{eqn:DistFlow_reconfig_SVI}, all constraints can be rewritten as linear equalities or inequalities (cf.~\eqref{eqn:bi2lin}). In contrast, in \eqref{eqn:DistFlow_reconfig_SVI} the product of continuous variables $\nu_k$ and $\ell_n$ cannot be equivalently replaced by any auxiliary variable as in \eqref{eqn:bi2lin} (the transformation requires at least one 0-1 binary variable in the product). Thus, problem~\eqref{opt:reconfig} is a mixed integer linear program with additional non-convex bilinear equality constraints.

Even though \eqref{opt:reconfig} resembles a classical model of the ROPF problem (if the objective is changed to loss minimization), we emphasize that the DG maximization objective in \eqref{opt:reconfig} leads to significant ramifications that render most of the typical solution methods ineffective. The details will be provided in Section~\ref{sec:pitfalls}. Instead, we propose to solve \eqref{opt:reconfig} using the SBB (i.e., spatial branch-and-bound) algorithm of Gurobi (version 9 or later), and we subsequently call this our ``proposed method''. The rationale of our choice is as follows:
\begin{enumerate}[label=(\alph*)]
	\item SBB returns (sub)optimal solutions with certifiable optimality gap.
	\item The SBB implementation in Gurobi is specialized in bilinear programming problem, allowing more effective use of the problem structure.
	\item Case studies with realistically sized benchmark instances indicate the practical computational performance of the proposed method.
\end{enumerate}
We note that different packages implementing the SBB algorithm (e.g., Gurobi, BARON, Knitro) behave differently for problem~\eqref{opt:reconfig}. The numerical experiment in Section~\ref{sec:SBB} indicates that Gurobi is the best performer among our choices.




\section{Illustration and Discussion of the Proposed Method} \label{sec:demo}
We evaluate the performance of the proposed method by solving instances of the ROPF problem in \eqref{opt:reconfig} derived from benchmarks listed in Table~\ref{tab:benchmarks}.
\begin{table}[h]
	\centering
	\caption{Distribution system benchmark data}
	\begin{tabular}{|c|c|c|c|c|c|}
		\hline
		name & \# of buses & \# of tie-switches & \# of DG & source\\
		\hline
		33bw & 33 & 5 & 2 & \cite{baran1989network,zimmerman2010matpower}\\
		\hline
		118zh & 118 & 15 & 10 & \cite{zimmerman2010matpower} \\
		\hline
		136ma & 136 & 21 & 7 & \cite{guimaraes2005reconfiguration,zimmerman2010matpower} \\
		\hline
		REDS 29+1 & 30 & 1 & 4 & \cite{ababei2010efficient} \\
		\hline
		REDS 83+11 & 84 & 13 & 4 & \cite{ababei2010efficient,su2003network} \\
		\hline
		REDS 135+1 & 136 & 21 & 4 & \cite{ababei2010efficient,guimaraes2005reconfiguration} \\
		\hline
		REDS 201+3 & 202 & 15 & 5 & \cite{guimaraes2005reconfiguration,ababei2010efficient} \\
		\hline
		REDS 873+7 & 874 & 27 & 5 & \cite{ababei2010efficient} \\
		\hline
	\end{tabular}
	\label{tab:benchmarks}
\end{table}

The topology, line parameters and load data of the benchmarks can be found in the respective sources in Table~\ref{tab:benchmarks}. For each benchmark we add DG sources in the system (location details omitted). The rating of each DG source is 10 MVA for all benchmarks except for REDS 873+7 where the rating is 50 MVA (large enough to prevent the capacity constraint \eqref{eqn:inverter_capacity} from being active). 

For each bus the active power generation lower and upper limits in \eqref{eqn:inverter_pg_bounds} are zero and the DG source apparent power rating, respectively. On the other hand, in \eqref{eqn:inverter_pg_bounds} reactive power generation is not sign-restricted, and its absolute value is no greater than the DG rating. For \eqref{eqn:inverter_pf}, the minimum power factor is 0.9 for all DG sources. Unless specified otherwise, for all buses the voltage lower and upper limits in \eqref{eqn:v_limits} are $0.95$ pu and $1.05$ pu, respectively. The current upper limit in \eqref{eqn:l_limit_heterogeneous}, on the other hand, varies from line to line. The values are in the order of a few hundred amperes. We also consider a variant of 33bw with line current rating enlarged to 600 A for all lines.

All computation in this and the next section is performed on a PC with an Intel i9-12900K CPU and 64 GB of RAM. The ROPF problem in \eqref{opt:reconfig} is solved with Gurobi 10.0.1 (with parameter ``non-convex'' set to 2) and we use AMPL to implement the optimization models. Gurobi is called with its default parameter values (e.g., target optimality gap 0.01\%).

\subsection{Reconfiguration with the 33-bus Benchmark}
\begin{figure}[h]
	\begin{center}
		\includegraphics[width=80mm]{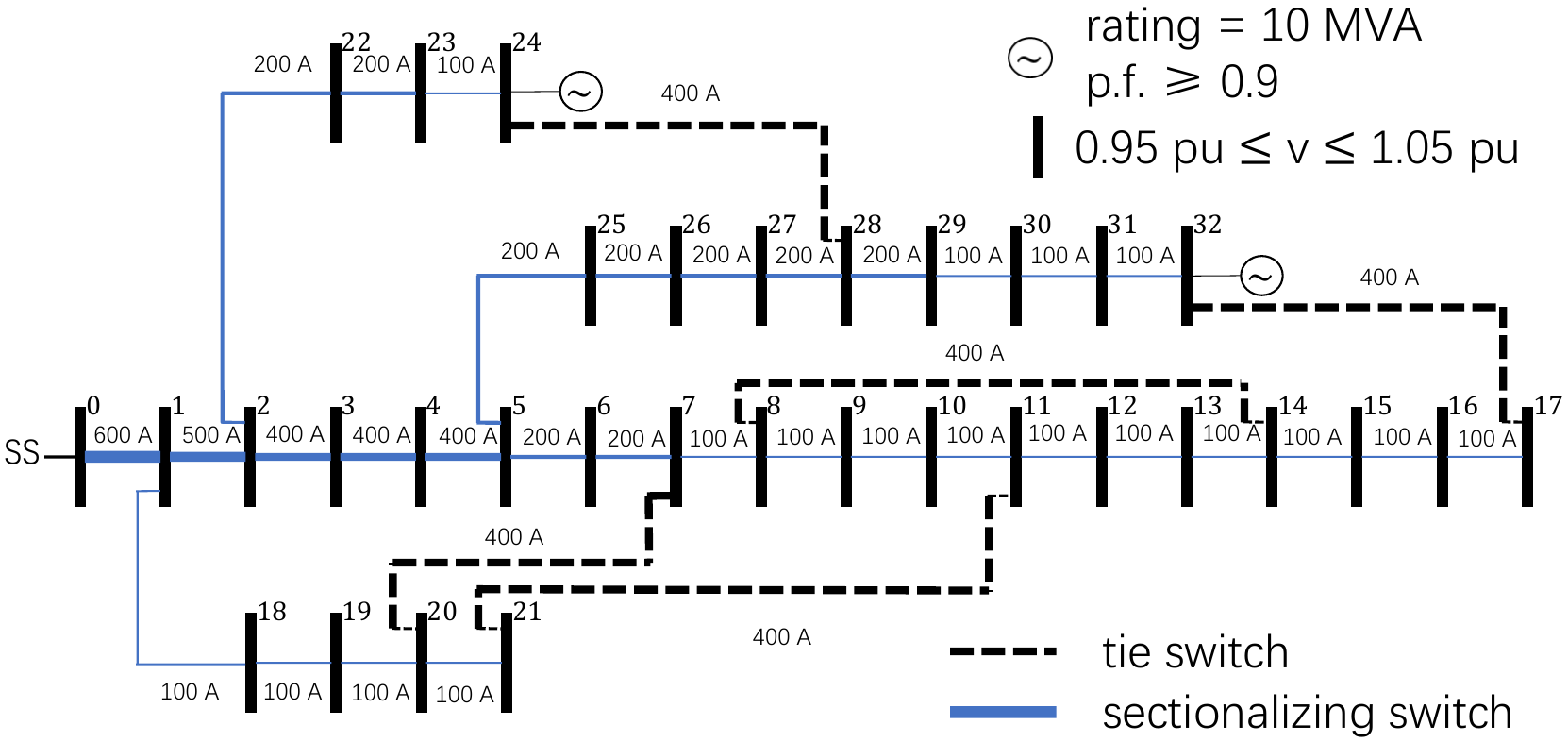}
		\caption{33-bus distribution system from \cite{baran1989network} with five tie-switches (dashed lines). Line width is proportional to current rating labeled next to the line.}
		\label{fig:33bw}
	\end{center}
\end{figure}

\begin{figure}[t]
	\begin{center}
		\includegraphics[width=80mm]{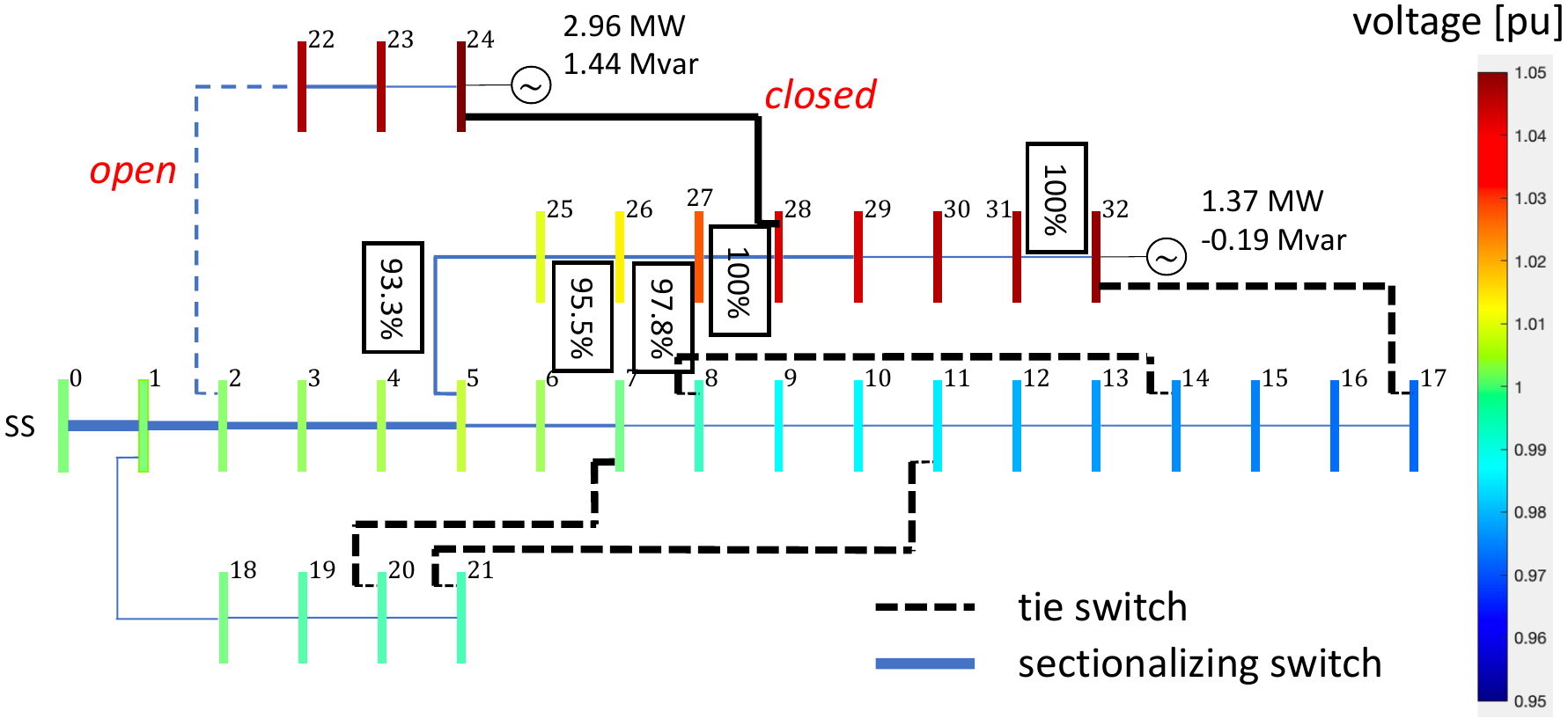}
		\caption{Reconfiguration with $K = 2$ (line \{2, 22\} opened and line \{24, 28\} closed). Bus color shows the voltage level and the percentage of current rating of bottleneck lines are shown. Total DG output is 4.33 MW.}
		\label{fig:33bw_K2}
	\end{center}
\end{figure}

\begin{figure}[t]
	\begin{center}
		\includegraphics[width=80mm]{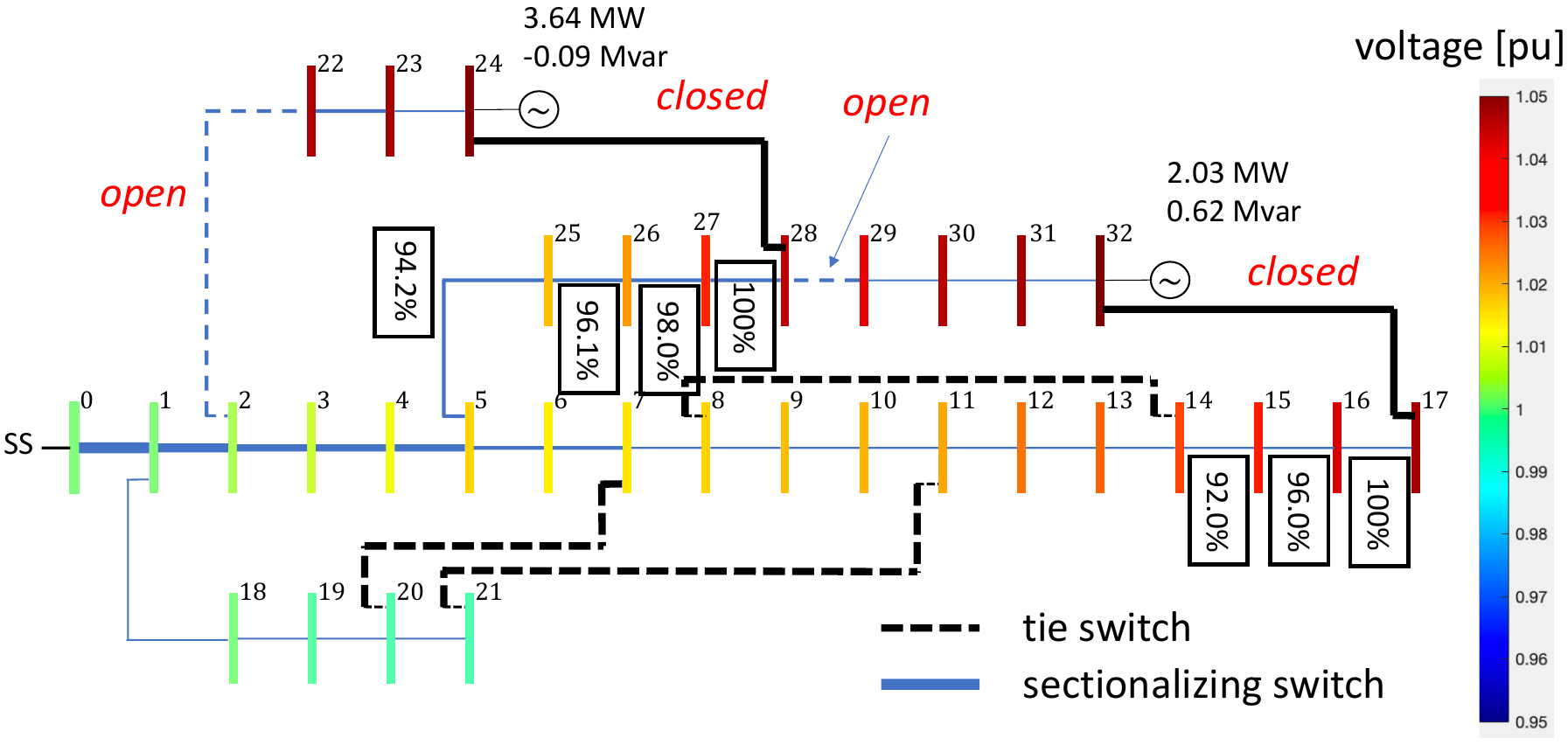}
		\caption{Reconfiguration with $K = 4$ (lines \{2, 22\}, \{28, 29\} opened, and lines \{24, 28\}, \{17, 32\} closed). Bus color shows the voltage level and the percentage of current rating of bottleneck lines are shown. Total DG output is 5.67 MW.}
		\label{fig:33bw_K4}
	\end{center}
\end{figure}

We consider the 33-bus benchmark 33bw described in Table~\ref{tab:benchmarks}. The topology, the locations of the two DG sources and the line current ratings are shown in Fig.~\ref{fig:33bw}.
We first attempt to solve the reconfiguration problem in \eqref{opt:reconfig} maximizing the total DG output with $K = 0$ (i.e., OPF with original network topology). The instance is infeasible due to operational constraints. However, if we relax the voltage lower limit from 0.95 pu to 0.9 pu, the instance becomes feasible and the voltages from bus 12 to bus 17 (towards the end of the feeder) are below 0.95 pu. 
Next, we solve \eqref{opt:reconfig} with $K = 2$ (i.e., one switch opened and another switch closed) with the usual voltage limits of 0.95 pu and 1.05 pu. The result is shown in Fig.~\ref{fig:33bw_K2}, with a total DG output of 4.33 MW.
In this reconfiguration, switch $\{2, 22\}$ is opened while switch $\{24, 28\}$ is closed. This raises the voltages from bus 5 to bus 17, satisfying the lower limit of 0.95 pu (e.g., the voltage at bus 17 is 0.972 pu). If we re-solve \eqref{opt:reconfig} with $K = 4$, the result is shown in Fig.~\ref{fig:33bw_K4}, with total DG output increased to 5.67 MW.
With $K = 4$, switch $\{28, 29\}$ is opened while switch \{17, 32\} is closed, in addition to the topology adjustments made in the case of $K = 2$. This prevents the two DG sources from sharing the path segment between bus 5 and bus 28, which is the bottleneck for $K = 2$. If we further increase $K = 6$ in \eqref{opt:reconfig}, the outcome is to open \{10, 11\} and close \{11, 21\}, cutting short the feeder from bus 0 to bus 17. This increases the total DG output to 5.72 MW. 



\subsection{Computation Time Evaluation}
The elapsed runtime to solve the benchmark instances is summarized in Table~\ref{tab:runtime}.
\begin{table}[!h]
	\centering
	\caption{Runtime [s] to achieve 0.01\% optimality gap, or optimality gap [\%] after one-hour runtime, for benchmark instances of \eqref{opt:reconfig}}
	\begin{tabular}{|c|c|c|c|c|c|c|}
		\hline
		name & $K = 0$ & $K = 2$ & $K = 4$ & $K = 6$ & $K = 8$\\
		\hline
		33bw & N/A & 1.56s & 1.69s & 2.21s & 4.06s \\
		\hline
		118zh & 0.82s & 20.18s & 28.7s & 105.5s & 275.74s \\
		\hline
		136ma & 0.14s & 4.88s & 12s & 15.35s & 58.16s \\
		\hline
		REDS 29+1 & 0.12s & 0.62s & 0.61s & 0.63s & 1.38s \\
		\hline
		REDS 83+11 & N/A & N/A & 0.42s & 2.33s & 3.32s \\
		\hline
		REDS 135+1 & 0.12s & 4.46s & 19.82s & 235.7s & 1.99\% \\
		\hline
		REDS 201+3 & 0.7s & 13.43s & 41.11s & 88.76s & 280.96s \\
		\hline
		REDS 873+7 & 5s & 265s & 1337s & 3.38\% & 4.54\% \\
		\hline
	\end{tabular}
	\label{tab:runtime}
\end{table}
In the table, N/A means that the instance is infeasible. In addition, for REDS 135+1 with $K = 8$ and REDS 873+7 with $K = 6$ and $K = 8$ the table entries show the percentage optimality gap after time limit of one hour. In general, runtime increases as the network size and $K$ increase. However, there are exceptions as REDS 135+1 appears to be difficult for $K = 8$ (cf.~136ma with $K = 8$ differing only in DG settings), whereas REDS 83+11 is easy (even with $K = 26$ the runtime is 11.8s, not shown in Table~\ref{tab:runtime}). Despite some irregularities in the runtime pattern, we conclude that the proposed method can reliably obtain guaranteed optimal DG set points and network configuration for moderately-sized instances (e.g., networks up to few hundred buses and $K \le 6$). However, discretion is needed when solving larger instances. Nevertheless, typically the marginal benefit diminishes rapidly with increasing $K$ (e.g., \cite{jabr2012minimum}).  Further, despite not being able to solve large instances to optimality, the proposed method can return suboptimal solutions with certified optimality gap. Runtime to achieve 5\% optimality gap for the benchmarks is summarized in Table~\ref{tab:runtime_5percent}.
\begin{table}[h]
	\centering
	\caption{Runtime to achieve 5\% optimality gap for benchmark instances of \eqref{opt:reconfig}}
	\begin{tabular}{|c|c|c|c|c|c|c|}
		\hline
		name & $K = 0$ & $K = 2$ & $K = 4$ & $K = 6$ & $K = 8$\\
		\hline
		33bw & N/A & 1.14s & 1.32s & 1.88s & 3.92s \\
		\hline
		118zh & 0.31s & 14.28s & 23.94s & 61.93s & 119.42s \\
		\hline
		136ma & 0.14s & 3.23s & 10.87s & 14.93s & 37.22s \\
		\hline
		REDS 29+1 & 0.11s & 0.61s & 0.61s & 0.63s & 0.63s \\
		\hline
		REDS 83+11 & N/A & N/A & 0.27s & 1.49s & 3.04s \\
		\hline
		REDS 135+1 & 0.11s & 3.09s & 8.80s & 5.96s & 91.74s \\
		\hline
		REDS 201+3 & 0.30s & 12.89s & 39.41s & 84.63s & 241.36s \\
		\hline
		REDS 873+7 & 0.58s & 86.87s & 482.4s & 1109s & 1854s \\
		\hline
	\end{tabular}
	\label{tab:runtime_5percent}
\end{table}
We conclude that, with reasonable target of optimality gap, the proposed method for DG maximizing ROPF is practical for networks of reasonable size and number of switch status changes. However, for very large number of switch status changes (or a problem to plan for possible line expansion), the proposed method should be augmented with more efficient methodologies.

\section{On the Choice of Optimization Solvers and Problem Formulations} \label{sec:pitfalls}
This section justifies the choice of our proposed method. We demonstrate that alternative methods for DG maximizing ROPF may not work as intended. These include interior point algorithm, MICP relaxation, linearized DistFlow model, alternative power flow equation formulations and SBB implementations.

\subsection{Pitfalls with Local Optimization Solvers} \label{subsec:local_opt}
We evaluate MATPOWER's OPF subroutine \texttt{runopf.m} to solve a special case of \eqref{opt:reconfig} with $K = 0$. To maximize the total DG output, MATPOWER is set to minimize $\sum_n c^0_n + c^1_n p_n^G + c^2_n (p_n^G)^2$ with $c^0_n = c^2_n = 0$ and $c^1_n = -1$ for all $n \in \mathcal{N}$. This is not the typical positive quadratic cost function. We run MATPOWER with the default options (e.g., flat start initial guess) and with two optimization solver choices: (a) MIPS interior point method and (b) Matlab's \texttt{fmincon} active set method. For the feasible instances in Table~\ref{tab:runtime}, MATPOWER solves the benchmark OPF instances with relative ease except that MIPS fails to converge for the 201+3 case.
However, MATPOWER results are sensitive to problem instance data and algorithmic choices. For example, for benchmark 33bw with line current ratings increased to 600 A for all lines and DG ratings decreased to 8 MW, the OPF instance becomes feasible but MATPOWER results vary greatly depending on the solver choice and initial guess strategy. Table~\ref{tab:matpower_OPF_33bw} shows the DG maximization results due to the following initial guess strategies \cite{zimmerman2010matpower}: 0 for flat start, 1 for ignoring the system state in the MATPOWER case, 2a and 2b for using the system state in the MATPOWER case (without and with update by first solving \eqref{opt:reconfig}) and 3a and 3b for solving power flow and using the resulting state (without and with update by first solving \eqref{opt:reconfig}). Table~\ref{tab:matpower_OPF_33bw} indicates that local optimization methods as in MATPOWER cannot reliably solve the DG maximization problem. {\R On the other hand, the proposed method returns the true optimal DG value of 10.45 MW.} 

The difficulty exhibited in Table~\ref{tab:matpower_OPF_33bw} is intrinsic to the DG maximization objective (as in \eqref{opt:reconfig}). If we instead minimize $\sum_n c^0_n + c^1_n p_n^G + c^2_n (p_n^G)^2$ with $c^0_n = c^2_n = 0$ and $c^1_n = 1$ (instead of $-1$) for all $n$, MATPOWER can reliably obtain the minimum total DG of 1.93 MW using both MIPS and \texttt{fmincon} for all choices of initial guess strategies.
\begin{table}[h]
	\centering
	\caption{33bw 600 A, MATPOWER OPF, max DG output}
	\begin{tabular}{|c|c|c|c|}
		\hline
		initial guess strategy & 0 & 1 & 2a \\
		\hline
		MIPS & 9.45 MW & 9.45 MW & 9.45 MW \\
		\hline
		\texttt{fmincon} & 10.44 MW & 10.44 MW & fail \\
		\hline
		\hline
		initial guess strategy & 2b & 3a & 3b \\
		\hline
		MIPS & 10.44 MW & 9.43 MW & 10.44 MW \\
		\hline
		\texttt{fmincon} & fail & 8.78 MW & fail \\
		\hline
	\end{tabular}
	\label{tab:matpower_OPF_33bw}
\end{table}
Therefore, when a ROPF problem with ``non-conventional'' objective (e.g., DG maximization) is considered, caution must be exercised regarding whether the optimization solver should work as intended.

\subsection{Pitfalls with MICP Relaxation} \label{subsec:socr}

\begin{figure}[!t]
	\begin{center}
		\includegraphics[width=50mm]{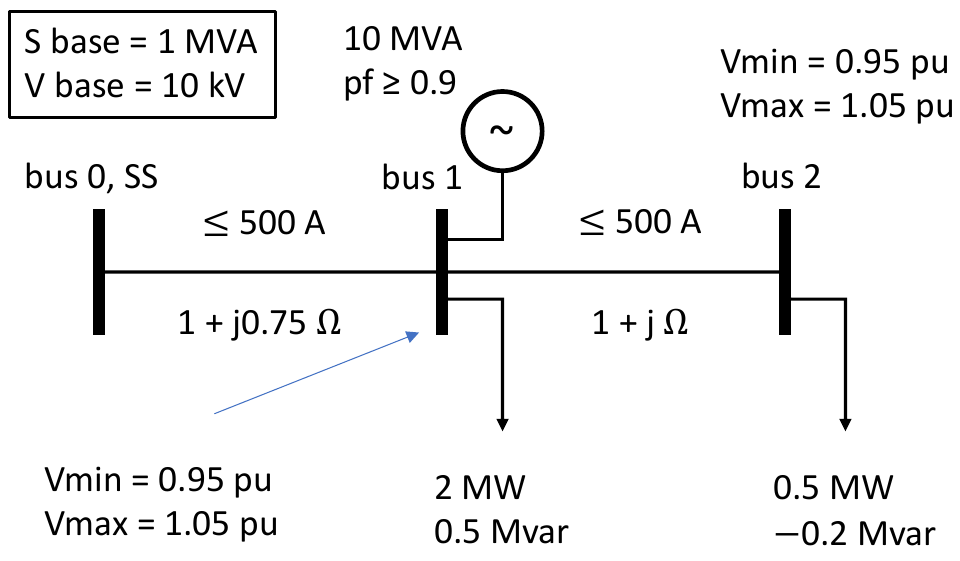}
		\caption{3-bus system showing failure of MICP relaxation. The voltage at bus 0 is 1 pu.}
		\label{fig:three_bus_ex}
	\end{center}
\end{figure}

An alternative to \eqref{opt:reconfig} is to replace the bilinear equations in \eqref{eqn:DistFlow_reconfig_SVI} with inequalities, leading to a relaxation \cite{farivar2013branch} as
\begin{subequations} \label{opt:SOCP_relax}
	\begin{align}
		\!\!\!\! \underset{\substack{\nu, \ell, p^G, q^G,  P, Q, \pi, \alpha}}{\text{maximize}} & \quad \sum\limits_{n \in \mathcal{N}} p^G_n \\
		\!\!\!\! \text{subject to} & \quad \text{\eqref{eqn:alpha_bounds} to \eqref{eqn:DistFlow_reconfig} except for \eqref{eqn:DistFlow_reconfig_SVI}, \eqref{eqn:inverter} to \eqref{eqn:l_limit_heterogeneous} and} \nonumber
		\\
		& \Big(\sum\limits_{k \in \bar{\mathcal{N}}} \pi_{kn} \nu_k \Big) \ell_n \ge P_n^2 + Q_n^2, \;\; n \in \mathcal{N} \label{eqn:DistFlow_reconfig_SVI_ineq}
	\end{align}
\end{subequations}
Unlike \eqref{eqn:DistFlow_reconfig_SVI}, inequalities \eqref{eqn:DistFlow_reconfig_SVI_ineq} allow \eqref{opt:SOCP_relax} to be categorized as a MICP problem amenable to efficient solution algorithms. However, a major question regarding \eqref{opt:SOCP_relax} is whether it is tight in the sense that \eqref{eqn:DistFlow_reconfig_SVI_ineq} hold as equations at optimality. If \eqref{opt:SOCP_relax} is not tight, the control decisions may violate grid constraints (e.g., voltage and current limits) as the DistFlow equations are not satisfied. For \eqref{opt:SOCP_relax} with ``standard'' objective such as total line loss minimization, the tightness of \eqref{opt:SOCP_relax} has been empirically and theoretically verified \cite{taylor2012convex,jabr2012minimum,kocuk2017new,farivar2013branch}. However, for DG maximization, the issue of tightness has not been explored in depth. In fact, our answer to this question is negative:   
we will present a three-bus counterexample to show that \eqref{opt:SOCP_relax} is not tight and the resulting control will lead to grid constraint violations (this also agrees with the observations in \cite{wei2017optimal}). 
Consider the OPF problem for the three-bus example in Fig.~\ref{fig:three_bus_ex}.
We use Gurobi to solve \eqref{opt:reconfig} and \eqref{opt:SOCP_relax} with $K = 0$. The results are summarized in Table~\ref{tab:three_bus_solutions}, {\R demonstrating that the proposed method is able to obtain the correct optimal DG outputs whereas MICP relaxation fails to do so (the different values of DG outputs are highlighted in red).}
\begin{table}[t]
	\centering
	\caption{Per unit OPF results for Fig.~\ref{fig:three_bus_ex} by \eqref{opt:reconfig} and \eqref{opt:SOCP_relax}}
	\label{tab:three_bus_solutions}
	\begin{tabular}{|c|c|c|c|c|c|c|}
		\hline
		\multicolumn{7}{|c|}{Bilinear program in \eqref{opt:reconfig}} \\
		\hline
		bus & $\nu$ & $\ell$ & $p^G$ & $q^G$ & $P$ & $Q$ \\
		\hline
		1 & 1.1025 & 25 & {\color{red} \bf 7.7518} & 0.39754 & $-4.9991$ & 0.092606 \\
		\hline
		2 & 1.0964 & 0.2645 & 0 & 0 & 0.50264 & $-0.19736$ \\
		\hline
		\hline
		\multicolumn{7}{|c|}{Conic relaxation in \eqref{opt:SOCP_relax}} \\
		\hline
		bus & $\nu$ & $\ell$ & $p^G$ & $q^G$ & $P$ & $Q$ \\
		\hline
		1 & 1.1025 & 25 & {\color{red} \bf 7.9991} & 0.64489 & $-4.9991$ & 0.092606 \\
		\hline
		2 & 1.0915 & 25 & 0 & 0 & 0.75 & 0.049999 \\
		\hline
	\end{tabular}
\end{table}
With the problem data in Fig.~\ref{fig:three_bus_ex} and the variable values in Table~\ref{tab:three_bus_solutions}, we can verify that the solutions of \eqref{opt:reconfig} and \eqref{opt:SOCP_relax} satisfy all constraints in the respective problems. For example, the reactive power balance equation in \eqref{eqn:DistFlow_Q} for bus 1 for \eqref{opt:reconfig} is $q^G_1 - q^L_1 = 0.39754 - 0.5 = -0.10246 = (-1) \times 0.092606 + (-0.19736) + 0.0075 \times 25 = -Q_1 + Q_2 + x_{01} \ell_1$.
The analogous equation holds for the solution of \eqref{opt:SOCP_relax} as well: $q^G_1 - q^L_1 = 0.64489 - 0.5 = 0.14489 = (-1) \times 0.092606 + 0.049999 + 0.0075 \times 25 = -Q_1 + Q_2 + x_{01} \ell_1$.
However, while equation \eqref{eqn:DistFlow_reconfig_SVI} at bus 2 holds for the solution of \eqref{opt:reconfig}: $\nu_1 \ell_2 = 1.1025 \times 0.2645 = 0.2916 = 0.50264^2 + {(-0.19736)}^2 = P_2^2 + Q_2^2$,
\eqref{eqn:DistFlow_reconfig_SVI_ineq} remains an inequality instead of equality for the solution of \eqref{opt:SOCP_relax}: $\nu_1 \ell_2 = 1.1025 \times 25 = 27.5625 > 0.565 = 0.75^2 + 0.049999^2 = P_2^2 + Q_2^2$.
Indeed, a posterior load flow analysis based on $(p^G, q^G)$ from \eqref{opt:SOCP_relax} reveals that the bus voltages are 1, 1.0539 and 1.0511 per unit, respectively. Similarly, the currents for the lines are 522.53 A and 51.235 A. The voltage and current upper limits are both violated by the solution of MICP relaxation. The aforementioned comparison is repeated for all benchmarks in Table~\ref{tab:benchmarks}. MICP relaxation solutions fail to satisfy equality \eqref{eqn:DistFlow_reconfig_SVI} in all cases, and result in voltage and current violations.

\begin{figure}[t]
	\begin{center}
		\includegraphics[width=50mm]{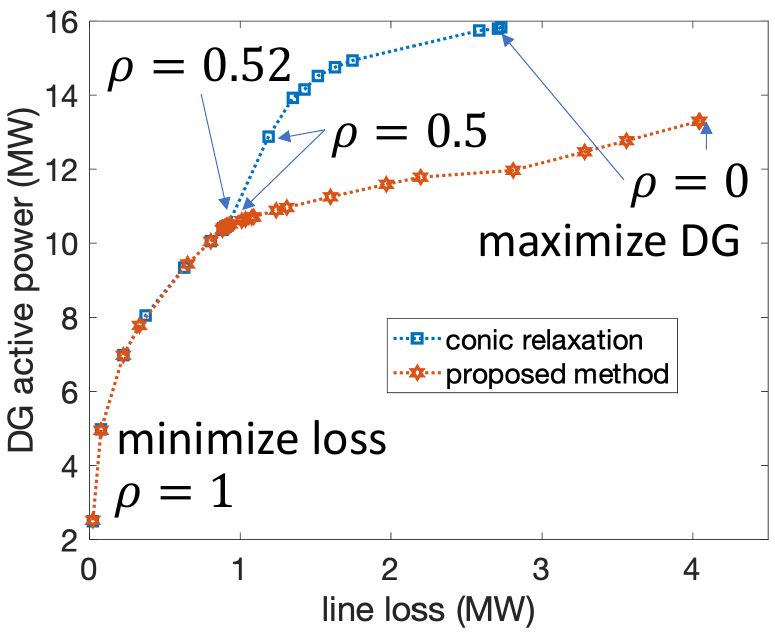}
		\caption{Pareto curves of 33bw benchmark (600A line rating) due to the proposed model in \eqref{opt:reconfig} and MICP relaxation in \eqref{opt:SOCP_relax}}
		\label{fig:33bw_pareto}
	\end{center}
\end{figure}

The objective functions in \eqref{opt:reconfig} and \eqref{opt:SOCP_relax} have a significant impact on the tightness of MICP relaxation. 
Following \cite{vcadjenovic2020maximization}, we consider the trade-off between DG maximization and line loss minimization by minimizing in \eqref{opt:reconfig} and \eqref{opt:SOCP_relax} the following weighted sum objective:
\begin{equation} \label{eqn:weighted_sum}
	\rho \sum\limits_{n \in \mathcal{N}} \big(\sum\limits_{k \in \bar{\mathcal{N}}} \pi_{kn} r_{kn}\big) \ell_n - (1 - \rho) \sum\limits_{n \in \mathcal{N}} p^G_n
\end{equation}
where $0 \le \rho \le 1$ is a trade-off parameter. 
With $\rho$ increasing from 0 to 1, the minimization of \eqref{eqn:weighted_sum} transitions from DG maximization to loss minimization. By solving \eqref{opt:reconfig} and \eqref{opt:SOCP_relax} with \eqref{eqn:weighted_sum} and running a posterior load flow analysis based on the resulting DG set points and network configuration, the trade-off can be visualized. Fig.~\ref{fig:33bw_pareto} shows the Pareto curves for the 33bw benchmark with universal line rating of 600 A and $K = 6$.
For $\rho > 0.5$, the two Pareto curves coincide since MICP relaxation solution satisfies \eqref{eqn:DistFlow_reconfig_SVI}. However, when $\rho \le 0.5$ MICP relaxation fails as it no longer depicts the true trade-off. 
This case study raises the awareness that the MICP relaxation approach should be used with caution especially when non-conventional objective function is considered.

\subsection{Pitfalls with Linearized DistFlow Approximation}
\begin{figure}[t]
	\begin{center}
		\includegraphics[width=50mm]{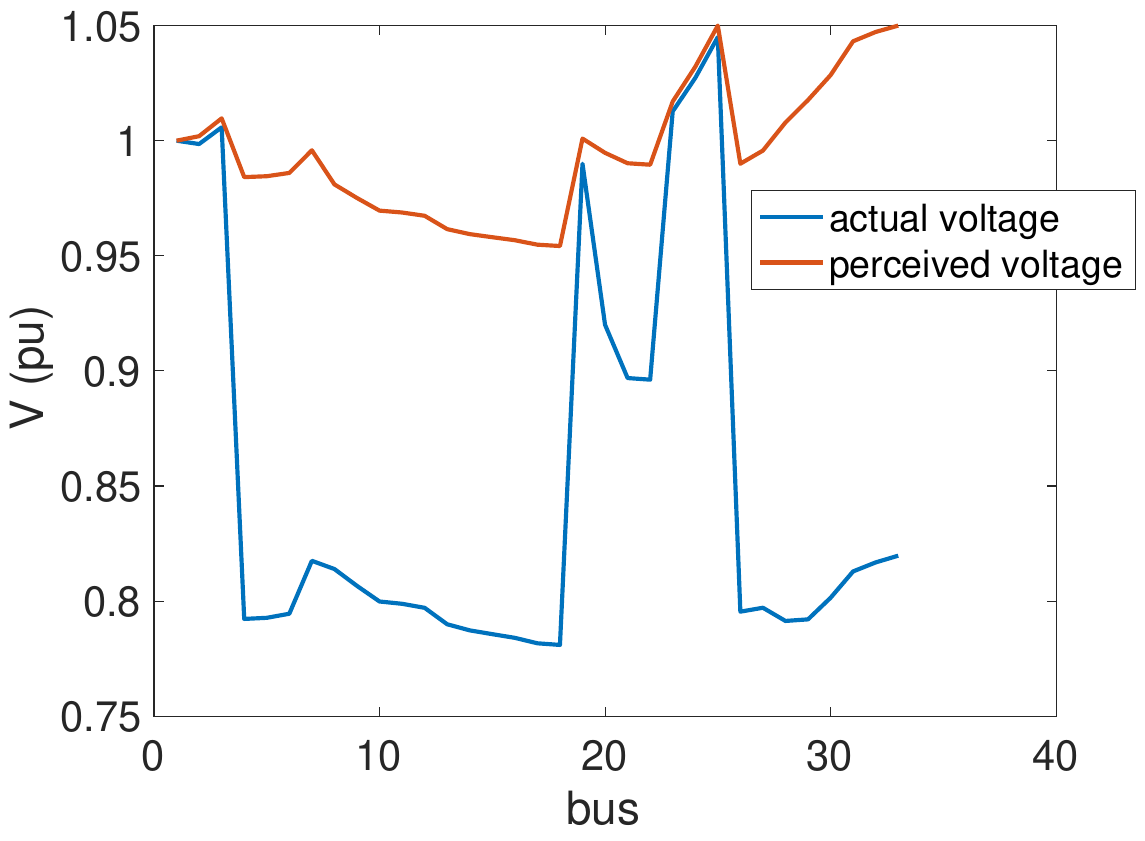}
		\caption{Benchmark 33bw (600 A line rating): actual voltage profile by a posterior load flow analysis and voltage profile perceived to be true by the LinDistFlow optimization model}
		\label{fig:LDF_33bw600A}
	\end{center}
\end{figure}

We evaluate a variant of \eqref{opt:reconfig} by replacing the DistFlow equations in \eqref{eqn:DistFlow_reconfig} with its linear approximation (i.e., \mbox{LinDistFlow} equations \cite{baran1989network}). The LinDistFlow variant ignores the loss-related terms multiplying $\ell_n$ in \eqref{eqn:DistFlow_reconfig_P}, \eqref{eqn:DistFlow_reconfig_Q} and \eqref{eqn:DistFlow_reconfig_Ohm}. Also, equation \eqref{eqn:DistFlow_reconfig_SVI} is dropped. Substituting the DistFlow equations with the LinDistFlow equations changes \eqref{opt:reconfig} to a MILP, which is relatively easy to solve. However, the voltage variables in the LinDistFlow variant overestimate the true voltages (e.g., \cite[eq.~(2a) and (3)]{SS_rollout2023}). Furthermore, since the squared current term $\ell$ is missing it is impossible to impose current upper limit constraints. Indeed, for our benchmark instances severe violation of current upper limit is experienced. To counter this, a surrogate current upper limit constraint can be imposed, similar to \eqref{eqn:DistFlow_reconfig_SVI}, as
\begin{equation} \label{eqn:LDF_lUB}
	\Big(\sum\limits_{k \in \bar{\mathcal{N}}} \pi_{kn} \nu_k \Big) \Big( \sum\limits_{k \in \bar{\mathcal{N}}} \pi_{kn} \overline{\ell}_{kn}\Big) \ge P_n^2 + Q_n^2, \quad n \in \mathcal{N}
\end{equation}
In \eqref{eqn:LDF_lUB} the product in the left-hand side can be equivalently replaced by linear terms as illustrated in \eqref{eqn:bi2lin}. Below we call LinDistFlow the variant of \eqref{opt:reconfig} with LinDistFlow equations and \eqref{eqn:LDF_lUB}. The inclusion of \eqref{eqn:LDF_lUB} relieves some undervoltage and overcurrent issues, but this is not enough. For instance, for the 33bw benchmark with 600 A line rating and $K = 2$, 
the actual voltages by a posterior load flow analysis and the voltages perceived to be true by the LinDistFlow optimization model (i.e., $(\nu)^{1/2}$) are very different as indicated by Fig.~\ref{fig:LDF_33bw600A}.
Also, current upper limits are violated in multiple lines (worst case violation being 115\% of rating) despite the surrogate constraint in \eqref{eqn:LDF_lUB}. This is likely to trigger an undervoltage trip. With reference to \eqref{eqn:DistFlow_reconfig_SVI}, overestimation of voltages will lead to underestimation of currents. Indeed, current upper limit are violated in multiple lines (worst case violation being 115\% of rating) despite the surrogate constraint in \eqref{eqn:LDF_lUB}. For LinDistFlow with $K = 4$ the result is even more extreme. The load flow analysis by Matpower fails to converge. {\R In contrast, solving \eqref{opt:reconfig} for the 33bw 600 A benchmark with $K = 2$ and $K = 4$, respectively, reliably leads to constraint-abiding control with guaranteed optimality.}

\subsection{Effect of Power Flow Equation Models} \label{sec:physics_models}
The model of power system physcis (e.g., DistFlow equations versus AC power flow equations) significantly affects the solver's performance to solve \eqref{opt:reconfig}. Here we investigate two ROPF models containing power flow equations other than the DistFlow equations \eqref{eqn:DistFlow_reconfig} adopted in our model in \eqref{opt:reconfig}: 
\begin{itemize}
	\item ACOPFC: ROPF model based on standard AC power flow equations by Capitanescu et.~al.~in \cite{capitanescu2014assessing}.
	\item ACOPFJ: ROPF model based on modified AC power flow equations by Jabr et.~al.~in \cite{jabr2012minimum}.
\end{itemize}
ACOPFC follows the model in \cite{capitanescu2014assessing}. On the other hand, ACOPFJ is based on the following power system physics model in \cite{jabr2012minimum}:
\begin{subequations} \label{eqn:Jabr_eqn}
	\begin{align}
		&p^G_n - p^L_n = \sum\limits_{k : (n,k) \in \mathcal{L}} P_{nk}, & &n \in \mathcal{N} \\
		&q^G_n - q^L_n = \sum\limits_{k : (n,k) \in \mathcal{L}} Q_{nk}, &  &n \in \mathcal{N} \\
		&P_{nk} = \sqrt{2} g_{nk} u^{(n,k)} - g_{nk} R_{nk} - b_{nk} T_{nk}, & &(n,k) \in \mathcal{L} \\
		&Q_{nk} = -\sqrt{2} b_{nk} u^{(n,k)} + b_{nk} R_{nk} - g_{nk} T_{nk}, & &(n,k) \in \mathcal{L} \\
		&2 u^{(n,k)} u^{(k,n)} \ge R_{nk}^2 + T_{nk}^2, & &(n,k) \in \mathcal{L} \label{eqn:Jabr_biineq} \\
		&R_{nk} = R_{kn}, & &(n,k) \in \mathcal{L} \\
		&T_{nk} = -T_{kn}, &  &(n,k) \in \mathcal{L} \\
		&0 \le u^{(n,k)} \le \frac{\overline{\nu}_n}{\sqrt{2}} \alpha_{nk}, &  &(n,k) \in \mathcal{L} \\
		&0 \le u_n - u^{(n,k)} \le \frac{\overline{\nu}_n}{\sqrt{2}} (1-\alpha_{nk}), & &(n,k) \in \mathcal{L}
	\end{align}
\end{subequations}
where $\mathcal{L}$ is the set of directed lines (normally open and normally closed) such that if $\{n,k\}$ is a line then $(n,k), (k,n) \in \mathcal{L}$, and $g_{nk} + j b_{nk}$ is the series admittance of line $(n,k)$. The decision variables include $P_{nk}$, $Q_{nk}$, $u^{(n,k)}$, $R_{nk}$, $T_{nk}$ and $u_n$ for all $(n,k) \in \mathcal{L}$. The ROPF formulation with \eqref{eqn:Jabr_eqn} imposed is a relaxation of \eqref{opt:reconfig} because of the inequalities in \eqref{eqn:Jabr_biineq}. However, replacing \eqref{eqn:Jabr_biineq} with equalities
\begin{equation} \label{eqn:Jabr_bieq}
	2 u^{(n,k)} u^{(k,n)} = R_{nk}^2 + T_{nk}^2, \quad (n,k) \in \mathcal{L}
\end{equation}
makes \eqref{eqn:Jabr_eqn} equivalent to the exact AC power flow equations (also equivalent to our DistFlow equations in \eqref{eqn:DistFlow_reconfig}). The ROPF model with \eqref{eqn:Jabr_eqn} (where the \eqref{eqn:Jabr_biineq} is replaced by \eqref{eqn:Jabr_bieq}) is the ACOPFJ formulation considered in our case study.

Both ACOPFC and ACOPFJ can be written as MINLP with bilinear equality constraints that is solvable using Gurobi. We test solving these models for the 33bw benchmark using Gurobi with time limit of 300 seconds. Table~\ref{tab:power_flow_model_33bw} shows the best DG output (if a feasible solution is found), the percentage optimality gap and the runtime to achieve the gaps. It is evident that for the DG maximizing ROPF, {\R our proposed model in \eqref{opt:reconfig} is the only practical choice when compared with ACOPFC and ACOPFJ.}
\begin{table}[h]
	\centering
	\caption{Comparing \eqref{opt:reconfig}, ACOPFC and ACOPFJ, 33bw}
	\begin{tabular}{|c|c|c|c|c|}
		\hline
		33bw & $K = 2$ & $K = 4$ & $K = 6$ & $K = 8$ \\
		\hline
		\multirow{2}{*}{\eqref{opt:reconfig}} & 4.34 MW & 5.67 MW & 5.72 MW & 5.72 MW\\
		\cline{2-5}
		& 0\%, 1.56s & 0\%, 1.69s & 0\%, 2.21s & 0\%, 4.06s \\
		\hline
		\multirow{2}{*}{ACOPFC} & N/A & N/A & N/A & 4.36 MW\\
		\cline{2-5}
		& fail, 300s & fail, 300s & fail, 300s & 359\%, 300s \\
		\hline
		\multirow{2}{*}{ACOPFJ} & N/A & 4.46 MW & 4.60 MW & 4.60 MW\\
		\cline{2-5}
		& fail & 174\%, 300s & 202\%, 300s & 191\%, 300s \\
		\hline
	\end{tabular}
	\label{tab:power_flow_model_33bw}
\end{table}

\subsection{Effect of SBB Solver Packages} \label{sec:SBB}
Different implementations of the SBB algorithm behave very differently when solving \eqref{opt:reconfig}. Here we evaluate three representative commercial packages:
	\begin{itemize}
			\item Gurobi 10.0.1: specialized for MINLP allowing only non-convex bilinear equations (and convex constraints).
			
			\item BARON 24.5.8: general SBB based global MINLP solver.
			
			\item Knitro 14.0.0: convex MINLP solver using SBB and interior point algorithms.
		\end{itemize}
	The three packages are run with default algorithmic parameter settings. Table~\ref{tab:SBB_33bw} shows the comparison result for the benchmark 33bw, showing how the three packages behave differently. The fastest runtime of Gurobi can be explained by the fact that it is specialized to consider only bilinear nonlinearity, whereas BARON is slowest as it aims to handle a much wider class of nonlinearities. Knitro is slower than Gurobi, and it may not return the true optimum as indicated by the MW results for $K = 2$ (4.14 MW versus 4.34 MW).
	\begin{table}[t]
			\centering
			\caption{Comparing Gurobi, BARON and Knitro, 33bw}
			\begin{tabular}{|c|c|c|c|c|}
					\hline
					33bw & $K = 2$ & $K = 4$ & $K = 6$ & $K = 8$ \\
					\hline
					\multirow{2}{*}{Gurobi} & 4.34 MW & 5.67 MW & 5.72 MW & 5.72 MW\\
					\cline{2-5}
					& 1.56s & 1.69s & 2.21s & 4.06s \\
					\hline
					\multirow{2}{*}{BARON} & 4.34 MW & 5.67 MW & 5.72 MW & 5.72 MW\\
					\cline{2-5}
					& 30.1s & 43.0s & 68.2s & 96.4s \\
					\hline
					\multirow{2}{*}{Knitro} & 4.14 MW & 5.67 MW & 5.71 MW & 5.72 MW\\
					\cline{2-5}
					& 8.64s & 8.68s & 16.5s & 49.5s \\
					\hline
				\end{tabular}
			\label{tab:SBB_33bw}
		\end{table}
	The result of a similar comparison case study with the 118zh benchmark is shown in Table~\ref{tab:SBB_118zh}, also indicating Gurobi's robustness for the purpose of solving \eqref{opt:reconfig}.
	\begin{table}[t]
			\centering
			\caption{Comparing Gurobi, BARON and Knitro, 118zh}
			\begin{tabular}{|c|c|c|c|c|}
					\hline
					118zh & $K = 0$ & $K = 2$ & $K = 4$ & $K = 6$\\
					\hline
					\multirow{2}{*}{Gurobi} & 27.47 MW & 28.57 MW & 29.44 MW & 30.17 MW\\
					\cline{2-5}
					& 0.82s & 20.2s & 28.7s & 105.5s \\
					\hline
					\multirow{2}{*}{BARON} & 27.47 MW & 28.56 MW & 29.27 MW & 30.15 MW\\
					\cline{2-5}
					& 2.10s & 570s & 2964s & 2hr (time limit) \\
					\hline
					\multirow{2}{*}{Knitro} & 27.46 MW & 28.18 MW & 29.35 MW & 29.49 MW\\
					\cline{2-5}
					& 0.41s & 183s & 1034s & 1176s \\
					\hline
				\end{tabular}
			\label{tab:SBB_118zh}
		\end{table}

\section{Practical Case Study - Real Network} \label{sec:real_network}

\subsection{533-bus Distribution System Description}
In this section, the ROPF problem \eqref{opt:reconfig} is applied to a real medium voltage distribution system in southern Sweden: the 533-bus network described in \cite{malmerthorin2023network}. With the expansion of DG, many DSOs wish to assess to what extent DG can be integrated in their systems without violating operational limits. This is referred to as performing a hosting capacity (HC) analysis. In \cite{ismael2019state-of-the-art}, four operational limits that constrain the HC in distribution networks are outlined -- voltage, overloading, protection and harmonics. With our proposed method we are imposing voltage and current limits, thereby dealing with the first two. Protection is partly included by requiring radial operation, but no further protection limits (i.e. short-circuit current) are set up. The system is assumed to be three-phase balanced without overtones, thus harmonics mitigation is not considered. Within these limits, the HC can be assessed and enhanced, using e.g. network reconfiguration, and our proposed ROPF method does this jointly.

\begin{figure}[h]
	\begin{center}
		\frame{\includegraphics[trim={4.5cm 1.8cm 3.3cm 1.3cm},clip,width=0.48\textwidth]{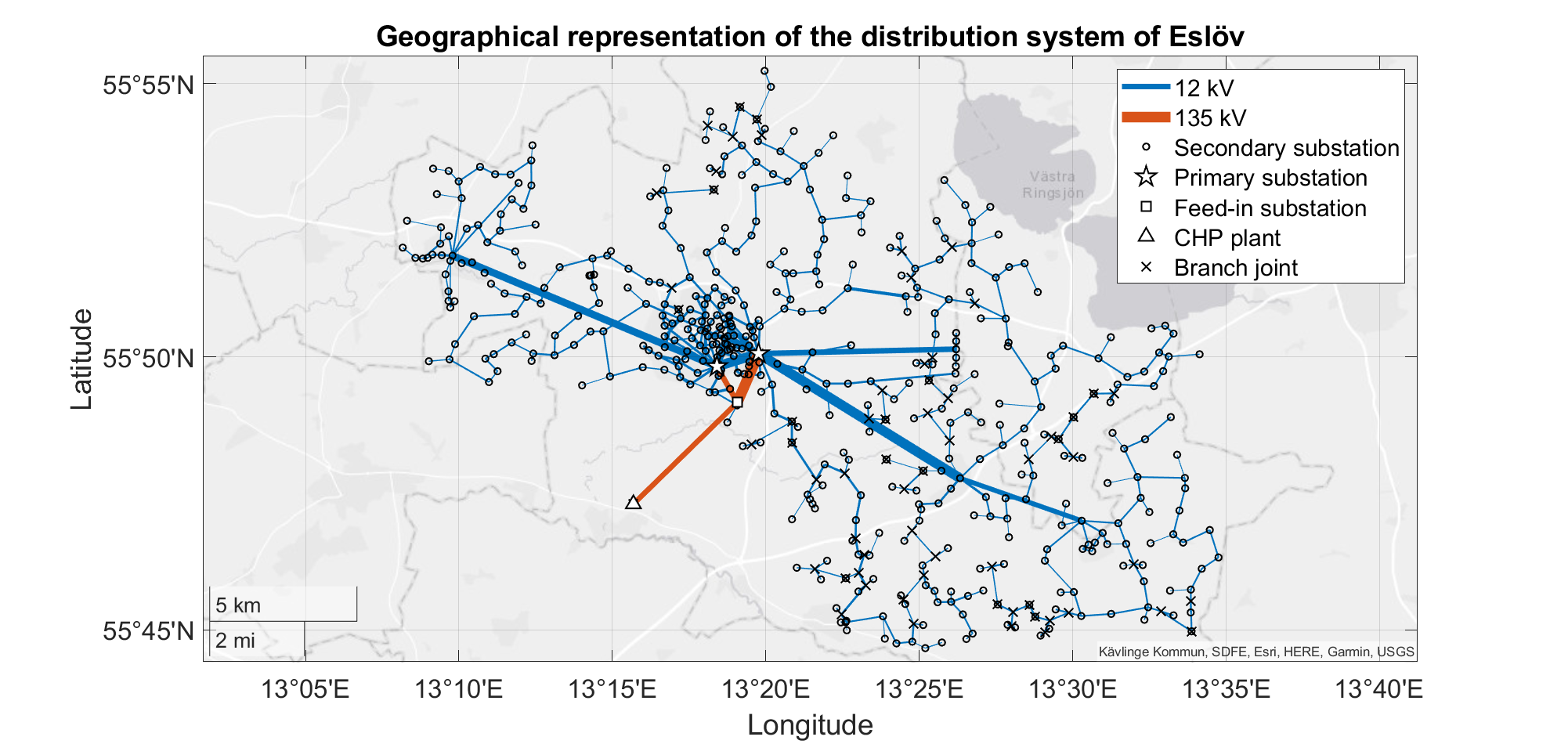}}
		\caption{Map of the 533-bus distribution system. The width of each line is proportional to its rated current.}
		\label{fig:533_bus}
	\end{center}
\end{figure}

The network in \cite{malmerthorin2023network} operates at 12 kV with a smaller part at 135 kV, covers about 20$\times$30 km and serves about 30,000 inhabitants plus an industrial area. It is connected to the regional 135 kV grid at a single feed-in station, modeled as a slack bus. A system overview is seen in Fig.~\ref{fig:533_bus}.

The system contains 533 nodes and 577 lines, thus radial operation \eqref{eqn:tree_num_of_edges} requires that 45 lines are always disconnected. Although a disconnected line might still be energized, the existence of an open point renders the entire line topologically switched off. The normal radial topology, used as base case ($K=0$) for subsequent simulations, can be seen in Fig.~\ref{fig:wind_scenario}. 

The existing wind and solar DG capacity amounts to 34.5 MW, installed at 262 nodes with capacities ranging from 3 kW to 2 MW per node. To avoid curtailment of already existing DG in the system, production from existing DG is embedded in $p^L_n$ by subtracting hourly generation from hourly load to obtain the net load data for each node. In the year 2022, the total net load in the system varied from 44.6 MW in Dec 16, 08:00 to -4.8 MW in Jul 10, 14:00. The latter will be referred to as the minimum net load hour.



The voltage constraints in \eqref{eqn:v_limits} are set to $\sqrt{\underline{\nu_{n}}} = 0.95$ pu and $\sqrt{\overline{\nu_{n}}} = 1.05$ pu, whereas the current constraints $\overline{\ell_{kn}}$ in \eqref{eqn:l_limit_heterogeneous} are set individually for each line, according to their respective rated currents. All system data was provided by the responsible DSO, \textit{Kraftringen}, which operates local distribution networks in multiple areas in southern Sweden, including the 533-bus network in question.
In the real network, there are often several line segments and conductor types connecting two substations, while in the model these are aggregated to single lines with equivalent rated currents and impedances. 
A full description of this aggregation procedure, and the model construction overall, can be found in \cite{malmerthorin2023network}. Once the aggregation has been made, the single lines have rated phase currents ranging from 93-1200 A, phase resistances from 0.003-1.86 $\Omega$ and phase reactances from 0.001-0.73 $\Omega$.

Furthermore, the real network is a three-phase system, whereas the optimization algorithm is designed for single-phase. Since the system is assumed to be three-phase balanced, this was managed by dividing all loads and generation capacities by a factor three, effectively treating each phase as a separate system. After simulations were run, the actual power flows were obtained by multiplying all power values post-optimization with a factor three. 
The final Matpower-model of the 533-bus system can be found at \texttt{\url{github.com/MATPOWER/matpower/blob/master/data/case533mt_lo.m}} with net load values from the minimum net load hour.

\subsection{Reconfiguration with 533-bus Example}
\label{subsec:533_results}
The ROPF problem in \eqref{opt:reconfig} is solved for the 533-bus system through two separate case studies. The first study involves the construction of a multiple node generation scenario, representing a potential future expansion of wind power. The second study is a spatially comprehensive DG expansion site analysis where additional DG capacity is added only to a single node in the system one at a time but iterating through the entire system, to see how high the HC is in each specific location and how much it can be increased through network reconfiguration.

In both the wind scenario and the DG expansion site analysis, the problem is solved with the minimum net load hour of 2022 as load input. The minimum net load hour is considered to be a good worst-case scenario, since bus voltages throughout the network are high during this hour and overvoltage is frequently the HC limiting factor. 

Additionally, in both the wind scenario and the DG expansion site analysis, results are given for $K=0$, $2$ and $4$ respectively. Naturally, the method can be applied also with higher $K$-values. However, this usually renders quickly diminishing marginal returns as shown in \cite{jabr2012minimum}, \cite{malmerthorin2023network} and Sec.~\ref{sec:windscenario}. In other words, the highest increases in HC are generally gained from the first switch events. Moreover, the optimization time rapidly increases as the search space grows exponentially with $K$. Finally, DSOs are in practice seldom interested in making too drastic topology changes in their networks based on specific load scenarios. The low $K$-scenarios reveal small topology changes from the normal configuration with potentially high rewards, which is attractive for grid planners.

\subsubsection{Wind scenario}

\label{sec:windscenario}
\begin{figure}[b]
	\begin{center}
		\frame{\includegraphics[trim={9.2cm 5cm 6.5cm 2.7cm},clip,width=0.45\textwidth]{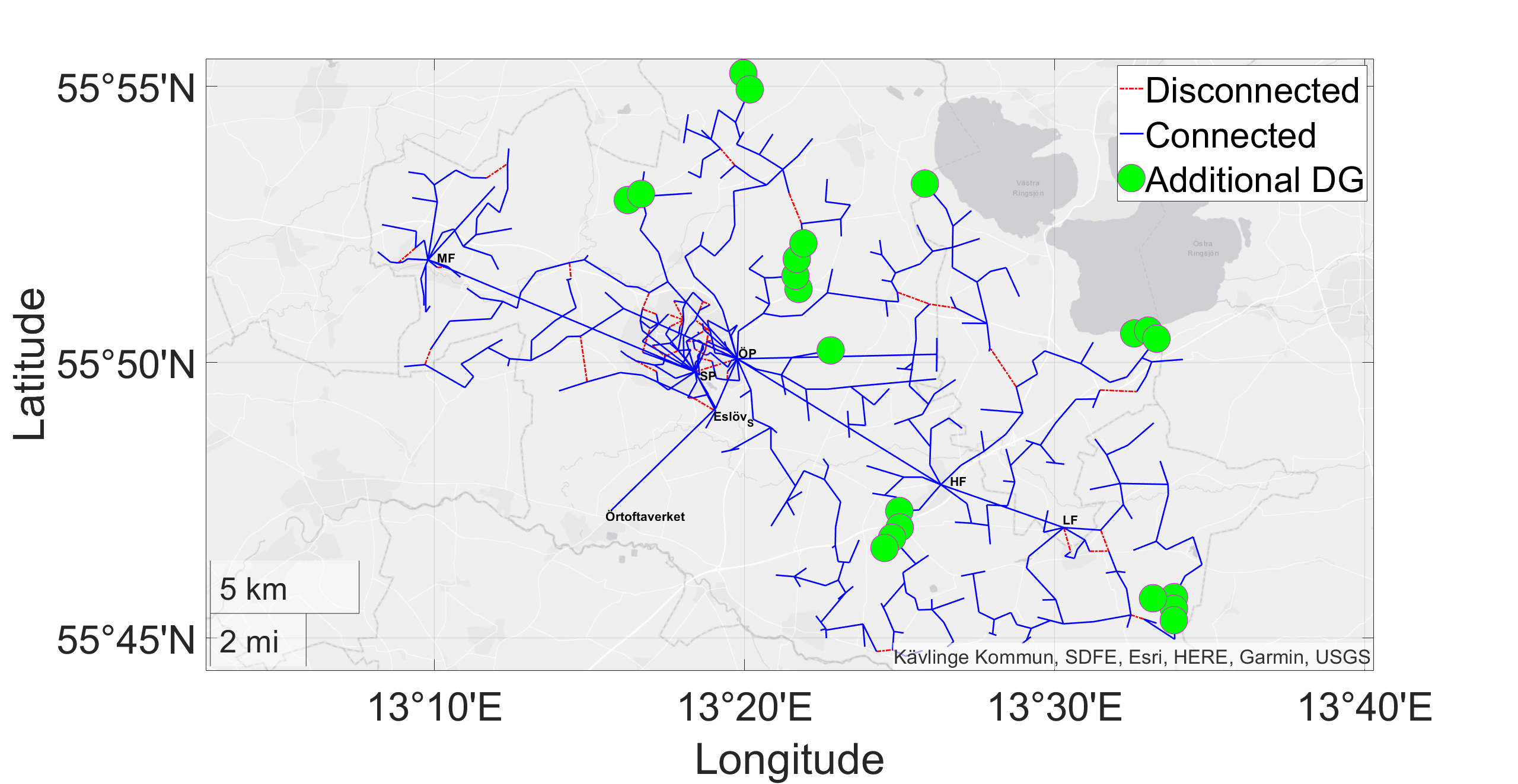}}
		\caption{Setup of the wind scenario in the normal radial topology. Green circles display the 21 nodes where additional DG units of 2 MW have been placed.}
		\label{fig:wind_scenario}
	\end{center}
\end{figure}

In the wind scenario, DG units of 2 MW are added to 21 rural nodes in clusters of up to four units. The placement of these clusters can be seen in Fig.~\ref{fig:wind_scenario}. The clusters are added without grid expansion or reinforcement, to determine the HC of the existing grid, before and after reconfiguration. Note that the 21 wind units are the only $p^G_n$ in the optimization, since existing DG is embedded in $p^L_n$. 

\begin{figure}[t]
	\begin{center}
		\includegraphics[width=120mm]{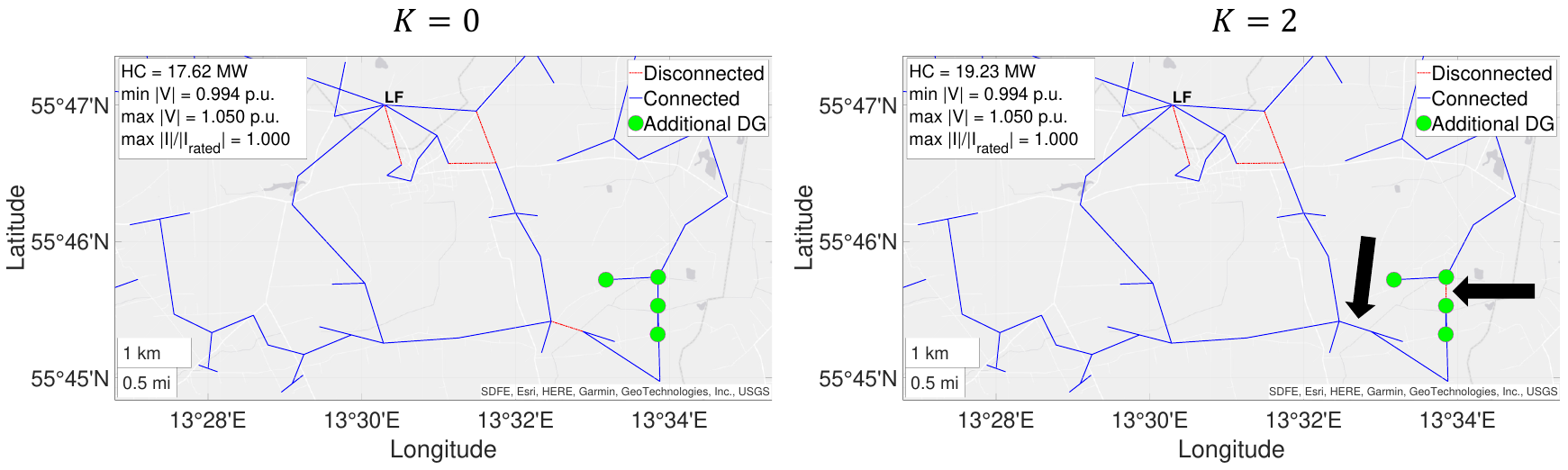}
		\caption{Left: Wind scenario in the normal radial topology ($K=0$), zoomed in at the southeastern part of the network. Right: Wind scenario after reconfiguration ($K=2$), the black arrows indicate where switch events occurred.}
		\label{fig:windscenario}
	\end{center}
\end{figure}

Solving \eqref{opt:reconfig} for this generation setup, $\sum_n p^G_n=17.6$, $19.2$ and $19.3$ MW for $K=0$, $2$ and $4$ respectively. As an initial observation, the added DG capacity of $21 \times 2=42$ MW is significantly curtailed both with and without reconfiguration. With $K=2$, the HC increases by 1.6 MW (or 9.1\%) compared to the normal topology. This increase is due to a pair of switch events that divides a cluster of four units in the southeast on two separate feeder lines, in contrast to the normal topology where they share a common feeder that turns out to be a bottleneck (see Fig.~\ref{fig:windscenario}). If four switch events are allowed ($K=4$), the algorithm returns the same first pair of switch events as with $K=2$, but also another pair of switch events closer to the distribution station LF. This has a more limited effect, increasing HC by only 0.1 MW from $K=2$.


\subsubsection{DG expansion site analysis}
In the DG expansion site analysis, additional DG capacity is added to a single node one at a time. This added capacity is large enough for the apparent power rating of the unit to not be an active constraint in \eqref{eqn:inverter_capacity}, allowing us to identify the HC of the grid. After adding DG to a single node, \eqref{opt:reconfig} is solved with $K=0$, $2$ and $4$. This procedure is repeated for all 532 non-slack nodes in the system, sometimes referred to as the iterative or detailed method \cite{stanfield2017optimizing}. 
No single optimization took more than 10 mins, which is considered reasonable for a site analysis.

\begin{figure}[h]
	\centering
	\includegraphics[trim={3.6cm 0.8cm 6cm 2.5cm},clip,width=0.49\textwidth]{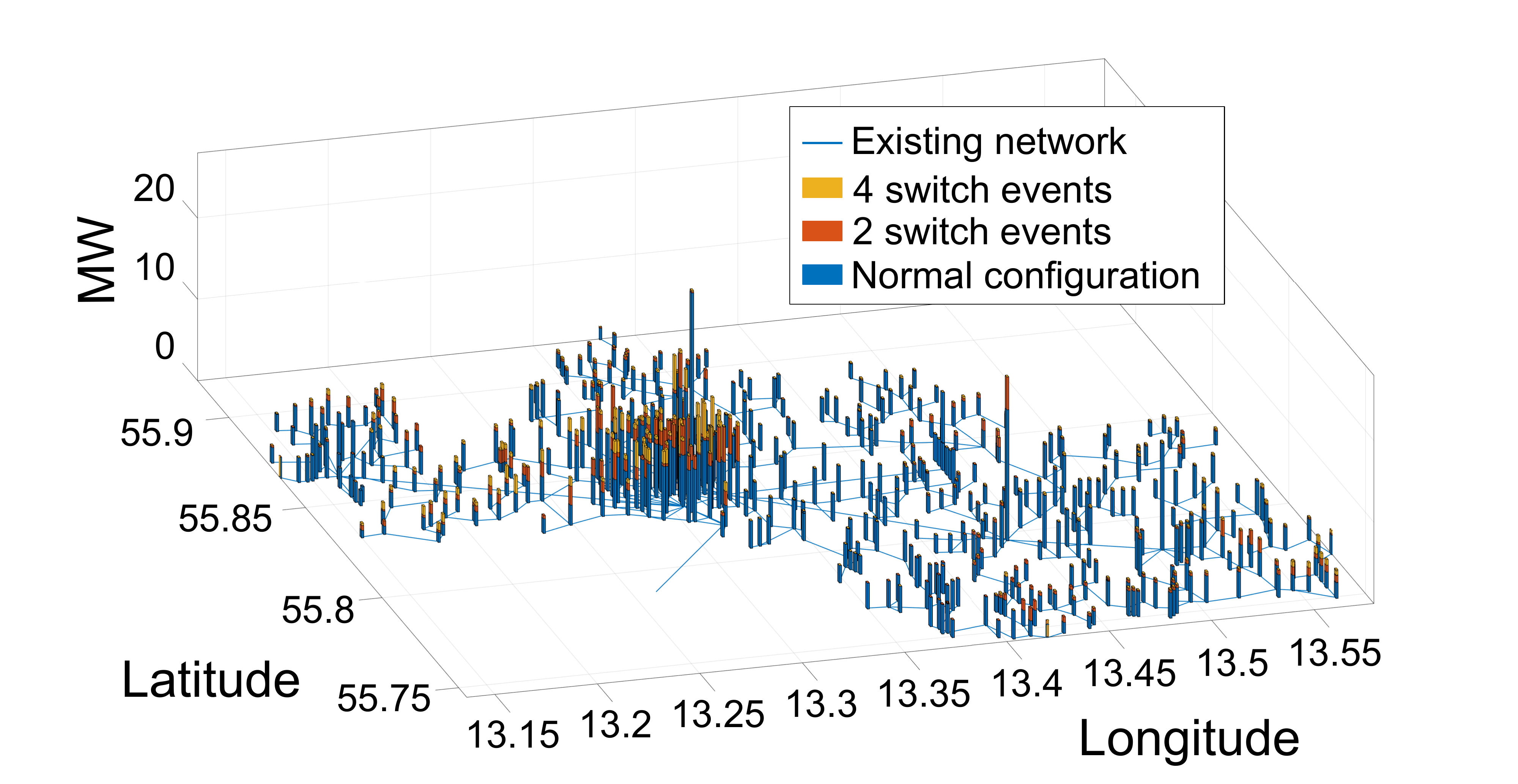}
	\caption{3D bar chart displaying the HC at each node. The blue, red and yellow parts are with $K=0$, $2$ and $4$ respectively.}
	\label{fig:3D_bar_chart}
\end{figure}

In Fig.~\ref{fig:3D_bar_chart}, the results from the $532 \times 3$ optimizations in the DG expansion site analysis can be seen.  The $x,y$-value indicates the location of a node in the network and the bar height is the HC at that node in MW. The blue part of each bar is the HC in the normal radial topology, whereas the red and yellow part is the additional HC including reconfiguration with $K=2$ and $4$ respectively. 

\begin{figure}[!t]
	\begin{center}
		\includegraphics[trim={2.7cm 0.4cm 6.3cm 2cm},clip,width=0.49\textwidth]{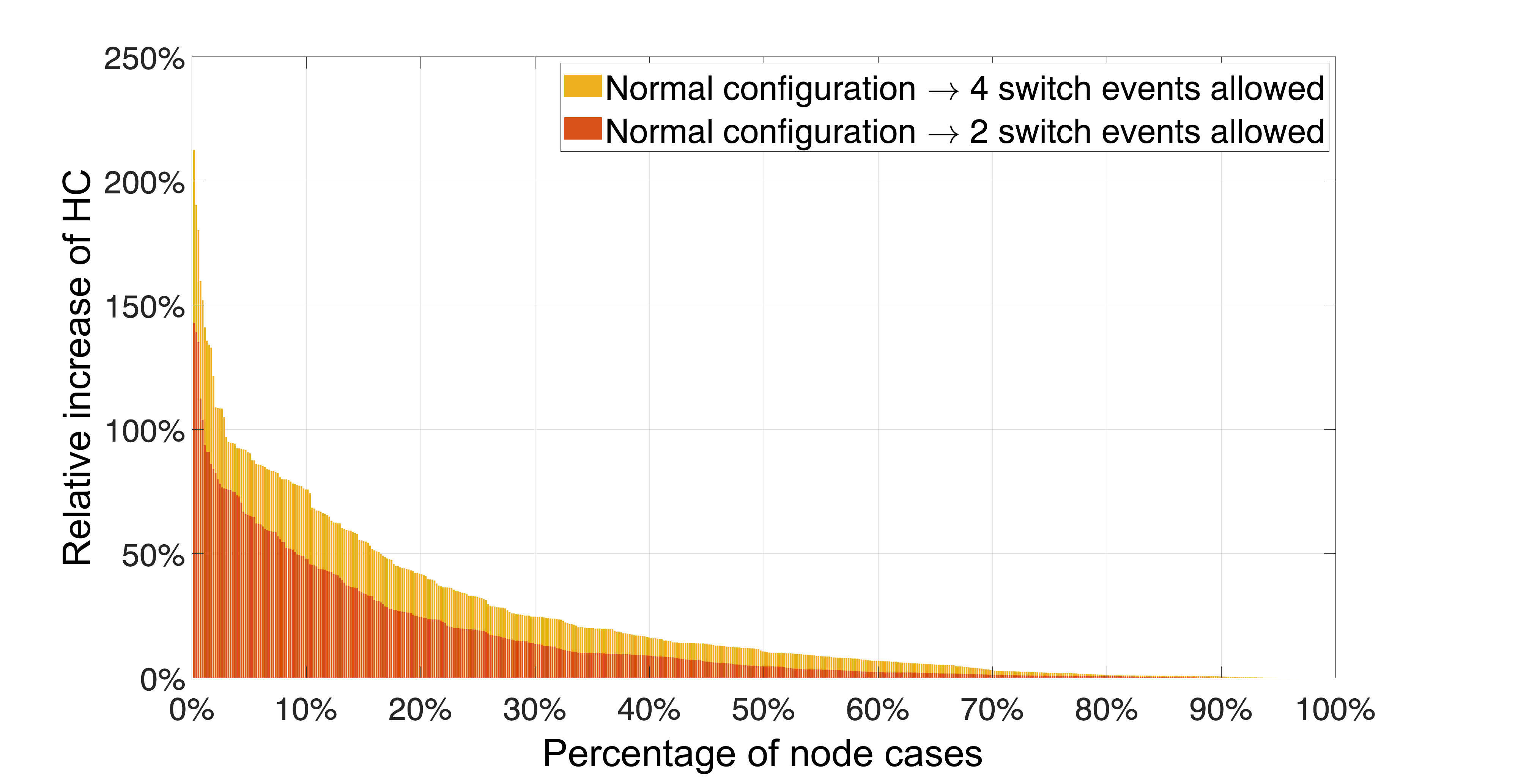}
		\caption{Share of nodes, $x\%$, that exhibit a HC increase of $y\%$ or higher from reconfiguration. Red depict the increase from $K=0$ to $K=2$, yellow from $K=0$ to $K=4$.}
		\label{fig:HC_increase}
	\end{center}
\end{figure}

Tall bars in Fig.~\ref{fig:3D_bar_chart} indicate a high HC, typically found at distribution substations or large industry nodes. On average, the HC is also higher at urban nodes than at rural. A large share of red or yellow in a bar indicates that the HC at that node can be improved substantially by reconfiguration. A large share of the urban nodes, but also some rural, exhibit this substantial increase in HC from reconfiguration. In Fig.~\ref{fig:HC_increase}, these relative HC improvements from network reconfiguration are presented explicitly. The figure displays how many of the nodes, $x\%$, that exhibit a relative HC increase of $y\%$ or higher from reconfiguration.

Fig.~\ref{fig:HC_increase} shows that the potential to enhance the HC through network reconfiguration varies considerably between nodes. The location of a node, and the network topology in its vicinity, has a significant effect on the benefit that can be achieved from reconfiguration. The greatest increase from $K=0$ to $K=4$ is as high as 212\% for a single node. Meanwhile, the median increase from $K=0$ to $K=4$ is only about 11\%.

%
%
\setlength{\textfloatsep}{10.0pt plus 2.0pt minus 2.0pt}
\section{Conclusion and Future Work} \label{sec:conclusion}
	DG maximization is more difficult than loss reduction for ROPF as the former objective renders standard techniques such as local optimization, MICP relaxation and linearization approximation ineffective, thus raising caution against their over-generalization and highlighting the importance of decision models featuring exact AC power flow equations (or equivalents). With a proper choice of problem formulation and optimization software package (e.g., DistFlow equations and Gurobi SBB), DG maximizing ROPF can indeed be handled reasonably well for moderately sized networks in a time scale relevant even for control center applications. 
	
	In this paper, the three-phase balanced problem has been considered. This is a natural starting point and a reasonable generalization as many distribution networks across the world are in balanced or close to balanced operation, the networks using three-phase feeders instead of single-phase laterals  \cite{mateo2018european},\cite{short2006electric}. For future work however, attention should be paid also to the unbalanced case, which is a reality in many networks and an important source of complexity. The ability of DG inverters and other sources of voltage control to compensate unbalances between phases could then be exploited further. Indeed, there exist variants of the DistFlow equations for three-phase unbalanced and meshed networks. Incorporating these equations in the ROPF model and identifying the proper solution approach are worth further investigation.
	
	

\bibliographystyle{elsarticle-num}
\bibliography{reconfig}

@article{baran1989network,
    title={Network reconfiguration in distribution systems for loss reduction and load balancing},
    author={Baran, Mesut E and Wu, Felix F},
    journal={IEEE Power Engineering Review},
    volume={9},
    number={4},
    pages={101--102},
    year={1989},
    publisher={IEEE}
}

@article{jabr2012minimum,
	title={Minimum loss network reconfiguration using mixed-integer convex programming},
	author={Jabr, Rabih A and Singh, Ravindra and Pal, Bikash C},
	journal={IEEE Transactions on Power systems},
	volume={27},
	number={2},
	pages={1106--1115},
	year={2012},
	publisher={IEEE}
}

@article{song2019new,
	title={A new formulation of distribution network reconfiguration for reducing the voltage volatility induced by distributed generation},
	author={Song, Yue and Zheng, Yu and Liu, Tao and Lei, Shunbo and Hill, David J},
	journal={IEEE Transactions on Power Systems},
	volume={35},
	number={1},
	pages={496--507},
	year={2019},
	publisher={IEEE}
}

@article{capitanescu2014assessing,
	title={Assessing the potential of network reconfiguration to improve distributed generation hosting capacity in active distribution systems},
	author={Capitanescu, Florin and Ochoa, Luis F and Margossian, Harag and Hatziargyriou, Nikos D},
	journal={IEEE Transactions on Power Systems},
	volume={30},
	number={1},
	pages={346--356},
	year={2014},
	publisher={IEEE}
}

@article{fu2018toward,
	title={Toward optimal multiperiod network reconfiguration for increasing the hosting capacity of distribution networks},
	author={Fu, Yang-Yang and Chiang, Hsiao-Dong},
	journal={IEEE Transactions on Power Delivery},
	volume={33},
	number={5},
	year={2018},
	publisher={IEEE}
}

@article{ababei2010efficient,
	title={Efficient network reconfiguration using minimum cost maximum flow-based branch exchanges and random walks-based loss estimations},
	author={Ababei, Cristinel and Kavasseri, Rajesh},
	journal={IEEE Transactions on Power Systems},
	volume={26},
	number={1},
	pages={30--37},
	year={2010},
	publisher={IEEE}
}

@article{zimmerman2010matpower,
	title={MATPOWER: Steady-state operations, planning, and analysis tools for power systems research and education},
	author={Zimmerman, Ray Daniel and Murillo-S{\'a}nchez, Carlos Edmundo and Thomas, Robert John},
	journal={IEEE Transactions on power systems},
	volume={26},
	number={1},
	pages={12--19},
	year={2010},
	publisher={IEEE}
}

@mastersthesis{malmerthorin2023network,
    title={Network reconfiguration for renewable generation maximization},
    author={Malmer, Gabriel and Thorin, Lovisa},
    year={2023},
    school={Lund University},
    type={Master's thesis}
}

@article{ismael2019state-of-the-art,
    title = {State-of-the-art of hosting capacity in modern power systems with distributed generation},
    journal = {Renewable Energy},
    volume = {130},
    pages = {1002-1020},
    year = {2019},
    issn = {0960-1481},
    doi = {https://doi.org/10.1016/j.renene.2018.07.008},
    author = {Sherif M. Ismael and Shady H.E. {Abdel Aleem} and Almoataz Y. Abdelaziz and Ahmed F. Zobaa},
}

@article{farivar2013branch,
	title={Branch Flow Model: Relaxations and Convexification -- \text{Part I}},
	author={Farivar, Masoud and Low, Steven H},
	journal={IEEE Transactions on Power Systems},
	volume={28},
	number={3},
	pages={2554--2564},
	year={2013},
	publisher={IEEE}
}

@techreport{stanfield2017optimizing,
    title={Optimizing the Grid: A Regulator’s Guide to Hosting Capacity Analyses for Distributed Energy Resources}, 
    author={Stanfield, Sky and Safdi, Stephanie and Mihaly, Shute and Weinberger, LLP}, 
    language={en},
    year={2017},
    institution={Interstate Renewable Energy Council}
}

@article{sou2022relaxed,
  title={Relaxed Connected Dominating Set Problem for Power System Cyber--Physical Security},
  author={Sou, Kin Cheong and Lu, Jie},
  journal={IEEE Transactions on Control of Network Systems},
  volume={9},
  number={4},
  pages={1780--1792},
  year={2022},
  publisher={IEEE}
}

@article{su2003network,
  title={Network reconfiguration of distribution systems using improved mixed-integer hybrid differential evolution},
  author={Su, Ching-Tzong and Lee, Chu-Sheng},
  journal={IEEE Transactions on power delivery},
  volume={18},
  number={3},
  pages={1022--1027},
  year={2003},
  publisher={IEEE}
}

@article{zhang2007improved,
  title={{An improved TS algorithm for loss-minimum reconfiguration in large-scale distribution systems}},
  author={Zhang, Dong and Fu, Zhengcai and Zhang, Liuchun},
  journal={Electric power systems research},
  volume={77},
  number={5-6},
  pages={685--694},
  year={2007},
  publisher={Elsevier}
}

@inproceedings{guimaraes2005reconfiguration,
  title={Reconfiguration of distribution systems for loss reduction using tabu search},
  author={Guimaraes, Marcos AN and Castro, Carlos A},
  booktitle={IEEE Power System Computation Conference (PSCC)},
  volume={1},
  pages={1--6},
  year={2005}
}

@article{taylor2012convex,
  title={Convex models of distribution system reconfiguration},
  author={Taylor, Joshua A and Hover, Franz S},
  journal={IEEE Transactions on Power Systems},
  volume={27},
  number={3},
  pages={1407--1413},
  year={2012},
  publisher={IEEE}
}

@article{home2022increasing,
  title={Increasing RES hosting capacity in distribution networks through closed-loop reconfiguration and Volt/VAr control},
  author={Home-Ortiz, Juan M and Macedo, Leonardo H and Vargas, Renzo and Romero, Rub{\'e}n and Mantovani, Jos{\'e} Roberto Sanches and Catal{\~a}o, Jo{\~a}o PS},
  journal={IEEE Transactions on Industry Applications},
  volume={58},
  number={4},
  pages={4424--4435},
  year={2022},
  publisher={IEEE}
}

@article{vcadjenovic2020maximization,
  title={Maximization of distribution network hosting capacity through optimal grid reconfiguration and distributed generation capacity allocation/control},
  author={{{C}}a{dj}enovi{c}, Rade and Jakus, Damir},
  journal={Energies},
  volume={13},
  number={20},
  pages={5315},
  year={2020},
  publisher={MDPI}
}

@article{SS_rollout2023,
    author = {Sou, Kin Cheong and Sandberg, Henrik},
    title = {Resilient Scheduling of Control Software Updates in Radial Power Distribution Systems},
    journal = {IEEE Transactions on Control of Network Systems},
    year = {2024},
    volume = {11},
    number = {3}
}

@article{low2014convex,
  title={Convex relaxation of optimal power flow—Part I: Formulations and equivalence},
  author={Low, Steven H},
  journal={IEEE Transactions on Control of Network Systems},
  volume={1},
  number={1},
  pages={15--27},
  year={2014},
  publisher={IEEE}
}

@inproceedings{sou2022joint,
  title={Joint Renewable Generation Maximization and Radial Distribution Network Reconfiguration},
  author={Sou, Kin Cheong and Gir{\'o}n, Kenny},
  booktitle={2022 IEEE PES Innovative Smart Grid Technologies-Asia (ISGT Asia)},
  pages={16--20},
  year = {2022}
}

@article{kocuk2017new,
  title={New formulation and strong MISOCP relaxations for AC optimal transmission switching problem},
  author={Kocuk, Burak and Dey, Santanu S and Sun, Xu Andy},
  journal={IEEE Transactions on Power Systems},
  volume={32},
  number={6},
  pages={4161--4170},
  year={2017},
  publisher={IEEE}
}

@article{wei2017optimal,
  title={Optimal power flow of radial networks and its variations: A sequential convex optimization approach},
  author={Wei, Wei and Wang, Jianhui and Li, Na and Mei, Shengwei},
  journal={IEEE Transactions on Smart Grid},
  volume={8},
  number={6},
  pages={2974--2987},
  year={2017},
  publisher={IEEE}
}

@inproceedings{hedman2011review,
  title={A review of transmission switching and network topology optimization},
  author={Hedman, Kory W and Oren, Shmuel S and O'Neill, Richard P},
  booktitle={2011 IEEE power and energy society general meeting},
  pages={1--7},
  year={2011},
  organization={IEEE}
}

@article{fisher2008optimal,
  title={Optimal transmission switching},
  author={Fisher, Emily B and O'Neill, Richard P and Ferris, Michael C},
  journal={IEEE Transactions on Power Systems},
  volume={23},
  number={3},
  pages={1346--1355},
  year={2008},
  publisher={IEEE}
}

@article{ergun2012transmission,
  title={Transmission system topology optimization for large-scale offshore wind integration},
  author={Ergun, Hakan and Van Hertem, Dirk and Belmans, Ronnie},
  journal={IEEE Transactions on Sustainable Energy},
  volume={3},
  number={4},
  pages={908--917},
  year={2012},
  publisher={IEEE}
}

@inproceedings{singh2022joint,
  title={Joint grid topology reconfiguration and design of watt-var curves for DERs},
  author={Singh, Manish K and Taheri, Sina and Kekatos, Vassilis and Schneider, Kevin P and Liu, Chen-Ching},
  booktitle={2022 IEEE Power \& Energy Society General Meeting (PESGM)},
  pages={1--5},
  year={2022},
  organization={IEEE}
}

@article{ramos2005path,
  title={Path-based distribution network modeling: application to reconfiguration for loss reduction},
  author={Ramos, E Romero and Exp{\'o}sito, A G{\'o}mez and Santos, J Riquelme and Iborra, F Llorens},
  journal={IEEE Transactions on power systems},
  volume={20},
  number={2},
  pages={556--564},
  year={2005},
  publisher={IEEE}
}

@article{lee2014robust,
  title={Robust distribution network reconfiguration},
  author={Lee, Changhyeok and Liu, Cong and Mehrotra, Sanjay and Bie, Zhaohong},
  journal={IEEE Transactions on Smart Grid},
  volume={6},
  number={2},
  pages={836--842},
  year={2014},
  publisher={IEEE}
}

@article{exposito1999reliable,
  title={Reliable load flow technique for radial distribution networks},
  author={Exp{\'o}sito, Antonio G{\'o}mez and Ramos, E Romero},
  journal={IEEE Transactions on Power Systems},
  volume={14},
  number={3},
  pages={1063--1069},
  year={1999},
  publisher={IEEE}
}

@article{jabr2006radial,
  title={Radial distribution load flow using conic programming},
  author={Jabr, Rabih A},
  journal={IEEE transactions on power systems},
  volume={21},
  number={3},
  pages={1458--1459},
  year={2006},
  publisher={IEEE}
}

@article{lavorato2011imposing,
  title={Imposing radiality constraints in distribution system optimization problems},
  author={Lavorato, Marina and Franco, John F and Rider, Marcos J and Romero, Rub{\'e}n},
  journal={IEEE Transactions on Power Systems},
  volume={27},
  number={1},
  pages={172--180},
  year={2011},
  publisher={IEEE}
}

@book{short2006electric,
publisher={CRC Press},
author={Thomas Allen Short},
title={Electric Power Distribution Equipment and Systems},
year={2005},
ISBN={9780367391676},
}

@article{mateo2018european,
	title={European representative electricity distribution networks},
	author={Mateo, Carlos and Prettico, Giuseppe and G{\'o}mez, Tom{\'a}s and Cossent, Rafael and Gangale, Flavia and Fr{\'\i}as, Pablo and Fulli, Gianluca},
	journal={International Journal of Electrical Power \& Energy Systems},
	volume={99},
	pages={273--280},
	year={2018},
	publisher={Elsevier}
}

@article{santos2022dynamic,
	title={Dynamic distribution system reconfiguration considering distributed renewable energy sources and energy storage systems},
	author={Santos, S{\'e}rgio F and Gough, Matthew and Fitiwi, Desta Z and Pogeira, Jos{\'e} and Shafie-khah, Miadreza and Catal{\~a}o, Jo{\~a}o PS},
	journal={IEEE Systems Journal},
	volume={16},
	number={3},
	pages={3723--3733},
	year={2022},
	publisher={IEEE}
}

@article{bahrami2024dynamic,
	title={Dynamic distribution network reconfiguration with generation and load uncertainty},
	author={Bahrami, Shahab and Chen, Yu Christine and Wong, Vincent WS},
	journal={IEEE Transactions on Smart Grid},
	volume={15},
	number={6},
	pages={5472--5484},
	year={2024},
	publisher={IEEE}
}
	
\end{document}